\documentclass{article}

\usepackage{PRIMEarxiv}
\usepackage[utf8]{inputenc}
\usepackage{textcomp}

\usepackage{graphicx}
\usepackage{amsmath}
\usepackage[version=4]{mhchem}
\usepackage{siunitx}
\usepackage{longtable,tabularx}
\usepackage{mathalpha}
\usepackage{bm}
\usepackage{adjustbox}
\usepackage{footnpag}			  
\usepackage{acronym}
\usepackage{caption}
\usepackage{subcaption}
\usepackage{siunitx}
\usepackage{booktabs}%
\usepackage{hyperref}
\usepackage{verbatim}
\usepackage{algorithm}
\usepackage{algpseudocode}
\usepackage{mathtools}
\usepackage[thinc]{esdiff}
\usepackage{capt-of}
\usepackage{longtable,tabularx}
\usepackage{caption}
\usepackage{subcaption}
\usepackage{arydshln}
\usepackage{setspace}
\usepackage{enumitem}
\usepackage{tikz}
\usetikzlibrary{shapes, arrows}
 \usetikzlibrary{calc} 
 \usetikzlibrary{fit}
 \usepackage[nocomma]{optidef}
 \usetikzlibrary{automata,positioning}
\usepackage{tablefootnote}

\newcommand{\pushcode}[1][1]{\hskip\dimexpr#1\algorithmicindent\relax}

\newcommand{\mat}{\begin{bmatrix}}
\newcommand{\matf}{\end{bmatrix}}
\newcounter{equationset}
\usepackage{mathrsfs}
\usepackage{amssymb}
\title{Convex Optimization-based Model Predictive Control for Active Space Debris 
Removal Mission Guidance
}

\author{
  Minduli C. Wijayatunga \\
  Phd Candidate\\
  Te Pūnaha Ātea - Space Institute \\
  The University of Auckland \\
  Auckland\\ \\
   \And
  Roberto Armellin\\
  Professor \\
  Te Pūnaha Ātea - Space Institute \\
  The University of Auckland \\
  Auckland\\ \\
     \And
  Laura Pirovano\\
  Research Fellow \\
  Te Pūnaha Ātea - Space Institute \\
  The University of Auckland \\
  Auckland\\\\
       \And
Harry Holt \\
  Research Fellow \\
  Te Pūnaha Ātea - Space Institute \\
  The University of Auckland \\
  Auckland\\\\
         \And
Claudio Bombardelli \\
 Professor \\
  Universidad Politécnica de Madrid \\ Spain \\}

\begin{document}
\maketitle

\begin{abstract}
A convex optimization-based model predictive control (MPC) algorithm for the guidance of active debris removal (ADR) missions is proposed in this work. A high-accuracy reference for the convex optimization is obtained through a split-Edelbaum approach that takes the effects of $J_2$, drag, and eclipses into account. When the spacecraft deviates significantly from the reference trajectory, a new reference is calculated through the same method to reach the target debris. When required, phasing is integrated into the transfer. During the mission, the phase of the spacecraft is adjusted to match that of the target debris at the end of the transfer by introducing intermediate waiting times. The robustness of the guidance scheme is tested in a high-fidelity dynamical model that includes thrust errors and misthrust events. The guidance algorithm performs well without requiring successive convex iterations. Monte-Carlo simulations are conducted to analyze the impact of these thrust uncertainties on the guidance. Simulation results show that the proposed convex-MPC approach can ensure that the spacecraft can reach its target despite significant uncertainties and long-duration misthrust events. 
\end{abstract}

\section*{Nomenclature}

{\renewcommand\arraystretch{1.0}
\noindent\begin{longtable*}{@{}l @{\quad=\quad} l@{}}
$T_{max}$ & Maximum thrust \\ 
$I_{sp}$&Specific impulse \\
$DC$ &  Duty Ratio \\ 
$m$ & spacecraft mass \\ 
$g_0$ & Gravitational acceleration (9.80665 m/s$^2$) \\
$J_2$ & $J_2$ constant ($1.08262668 \times 10^{-3}$)\\
$\mu$& Gravitational parameter ($\SI{3.986e5}{\kilo\meter\cubed\per\second\squared}$)\\  
$t$ & time \\ 
$\Delta v'$ & Gauss variational equations-based estimate of the $\Delta v$ required to make an orbital change \\ 
$\Delta v'_{fp}$ & $\Delta v'$ obtained from the forward propagation of the convex solution under realistic conditions. \\
$\mathscr{a}$ & Acceleration \\  
$\eta$ & thrust profile \\
$N_{run}$ & index of $t^{convex}$ \\ 
$N_{seg}$ & length of a tracking segment \\

\multicolumn{2}{@{}l}{Subscripts}\\
$0$ & initial state \\
$T$ & target state\\
$f$ & reached state during tracking \\ 
$d$	& drift orbit parameters from the RAAN matching scheme \\ 
$convex$ & parameters related to the convex tracking \\ 
$PMDT$ & parameters related to the preliminary mission design tool (PMDT)\\ 
$fp$ & parameters related to the forward propagation after convex optimization \\
$req$ & real spacecraft engine specifications \\
$CC$ & Cartesian coordinates \\ 
$COE$ & Classical orbital elements/ Keplerian coordinates\\
$MEE$ & Modified equinoctial elements \\
$GEq$ & Generalized equinoctial elements \\ 
\end{longtable*}}

\section{Introduction}

The low Earth orbit (LEO) space environment is becoming increasingly congested with space debris. As of September 2023,  rocket or payload debris makes up  {43\% of the LEO population} \cite{esaspacedebrisofficeESAANNUALSPACE2023}, and the average rate of debris collisions has increased to four or five objects per year \cite{Maestrini2021}. As satellites become increasingly essential to daily life, more are added to expand space-enabled services. However, additional launches increase the risk of collision for all satellites as they further saturate space with objects, endangering critical space infrastructure such as the International Space Station in LEO. Recent studies have shown that ensuring adequate post-mission disposal of new satellites is insufficient to prevent a future collision cascade and that active debris reduction by five to ten large objects per year is necessary \cite{CASTRONUOVO2011848}.  \par 
Active debris removal (ADR) is defined as removing derelict objects from space, thus minimizing the build-up of unnecessary objects and lowering the probability of on-orbit collisions  \cite{bonnal2013active,liou2010parametric}. ADR has become significant over the past two decades, leading to numerous studies and implementations of potential debris removal missions and technologies. The End of Life Service by Astroscale-demonstration (ELSA-d) mission was launched in March 2021 and has successfully tested both rendezvous algorithms required for ADR and a magnetic capture mechanism to remove objects carrying a dedicated docking plate at the end of their missions  \footnote{\url{https://astroscale.com/elsa-d-mission-update/}}. The RemoveDebris mission by the University of Surrey is another project that demonstrated various debris removal methods, including harpoon and net capture \cite{Forshaw2020TheAS}.\par 

While these individual removal missions are essential milestones towards ADR implementation, large-scale missions that target multiple objects might be necessary to compete with the current debris growth rate \cite{liou2010parametric, white2014}. Such missions are expected to rely significantly on autonomy to reduce costs and maintain reliability. However, enabling technologies for such missions are still under development \cite{surveydebris}.\par 

In our previous work \cite{ADRMW}, a preliminary trajectory optimization tool (PMDT) was developed to obtain fuel and time-optimal trajectories for missions that remove multiple debris from orbit. The PMDT considers the impact of the oblateness of the Earth ($J_2$), eclipses-similar to \cite{4}-, duty cycle   and drag, and can optimize a three-debris removal mission in under a minute. Firstly, it utilizes a version of the classical Edelbaum's method \cite{2} to calculate the optimal time of flight and fuel expenditures of a single transfer. Then, a right ascension of the ascending node (RAAN) matching algorithm introduces an intermediate drift orbit where the spacecraft can utilize  $J_2$ perturbations to reach a desired RAAN. The mission's fuel consumption and flight time are optimized by adjusting the drift orbit and the launch time. In \cite{ADRMW}, PMDT results were used as a reference to provide guidance for complex multi-ADR mission profiles via existing guidance laws such as the $\Delta v$ law\cite{slim} and Q-law \cite{Petropoulos2005}. This showed that these guidance laws could track the PMDT reference relatively well and that the $\Delta v$  and time of flight ($TOF$) obtained from the PMDT are reasonable estimates for the transfers. However, the classical guidance laws required slightly more fuel than estimated by the PMDT and could only achieve an accuracy of $\sim$10 km in the semi-major axis and $\sim$0.1 deg in inclination. The limited accuracy of the classical guidance laws was speculated to be due to their heavy reliance on approximations and simplifications (i.e., the maximum rate of change approximation of the Gauss variational equations used in the Q-law \cite{Petropoulos2005} and the simplified $\Delta v$ estimation formulae used in the $\Delta v$-law \cite{slim}). \par

In contrast, Model Predictive Control (MPC) can provide significantly more accurate control, as it can account for perturbations in real-time \cite{Bashnick}. It can also adapt for significant divergences, as it takes a receding horizon control approach \cite{MPC1}. The use of MPC in the context of spacecraft rendezvous has been recently explored by L. Ravikumar \cite{RAVIKUMAR2020518}, C. Bashnick \cite{Bashnick} and R. Vazquez \cite{VAZQUEZ2015251}, who have all identified solving optimization problems at each of the control intervals as a computationally complex and time-consuming process. Convex optimization is appealing in this context as a single iteration of convex optimization can obtain optimal solutions in polynomial time \cite{hu2020convex,ack}. The use of convex optimization within MPC in spacecraft trajectory optimization has also been explored frequently in the literature \cite{doi:10.2514/1.G000218, WANG2023477,sun2018convex}. However, it is noted that convex optimization relies on the convexification of nonconvex dynamics and constraints, leading to the need for successive convexifications. {In the presence of highly nonconvex dynamics, this brings about inefficient optimizations that do not reach the global optimal solution \cite{hu2020convex}.} \par 

This study focuses on providing autonomous guidance for multi-ADR missions using a novel convex-based MPC method. A reference is first calculated using the PMDT discussed in \cite{ADRMW}, which is loosely tracked by a convex-based optimization in predefined time segments. Intermediate tracking constraints of the convex optimization are formulated as soft constraints to avoid unnecessary fuel consumption. If the spacecraft is expected to rendezvous with a target, it needs to match the phase angle of that target. In this work, phasing is embedded in the transfer itself by simply shifting the endpoint of the convex segments by a calculated amount of time. This allows the spacecraft to reach a targetted phase angle without consuming any fuel for phase-matching. {Lastly, if the real trajectory and the reference start to deviate from each other due to the accumulation of errors over time, a new reference is recomputed from the current position to the target using the PMDT in an MPC manner. }  \par 

{Note that a convexified trajectory can differ from reality due to convexification errors and discrepancies between reality and the dynamical model used. Performing successive convexifications mitigates the former but does not affect the latter. Furthermore, successive convexifications are unideal for autonomy as they may be time-consuming, given the uncertainty regarding the number of iterations required for convergence. If the initial convexification is performed about a highly accurate reference, it can enhance the region of validity of the convex formulation \cite{boyd_vandenberghe_2004}, reducing the need for successive iterations. In this work, such a reference is generated via the PMDT. The capacity for trajectory recomputation also helps resolve modeling and convexification error build-up without resorting to successive iterations. }   \par 

The remainder of this paper is organized as follows. Section 2 outlines the overview of the multi-ADR mission for which the guidance is applied. The proposed MPC guidance is described in depth in Section 3. Section 4 discusses the results of a fuel-optimal tour to remove a discarded rocket body, and Section 5 provides the concluding remarks.

\section{Mission Overview}

\begin{figure}[hbt!]
\begin{minipage}[b]{0.6\linewidth}
    \centering
    \includegraphics[width = \textwidth]{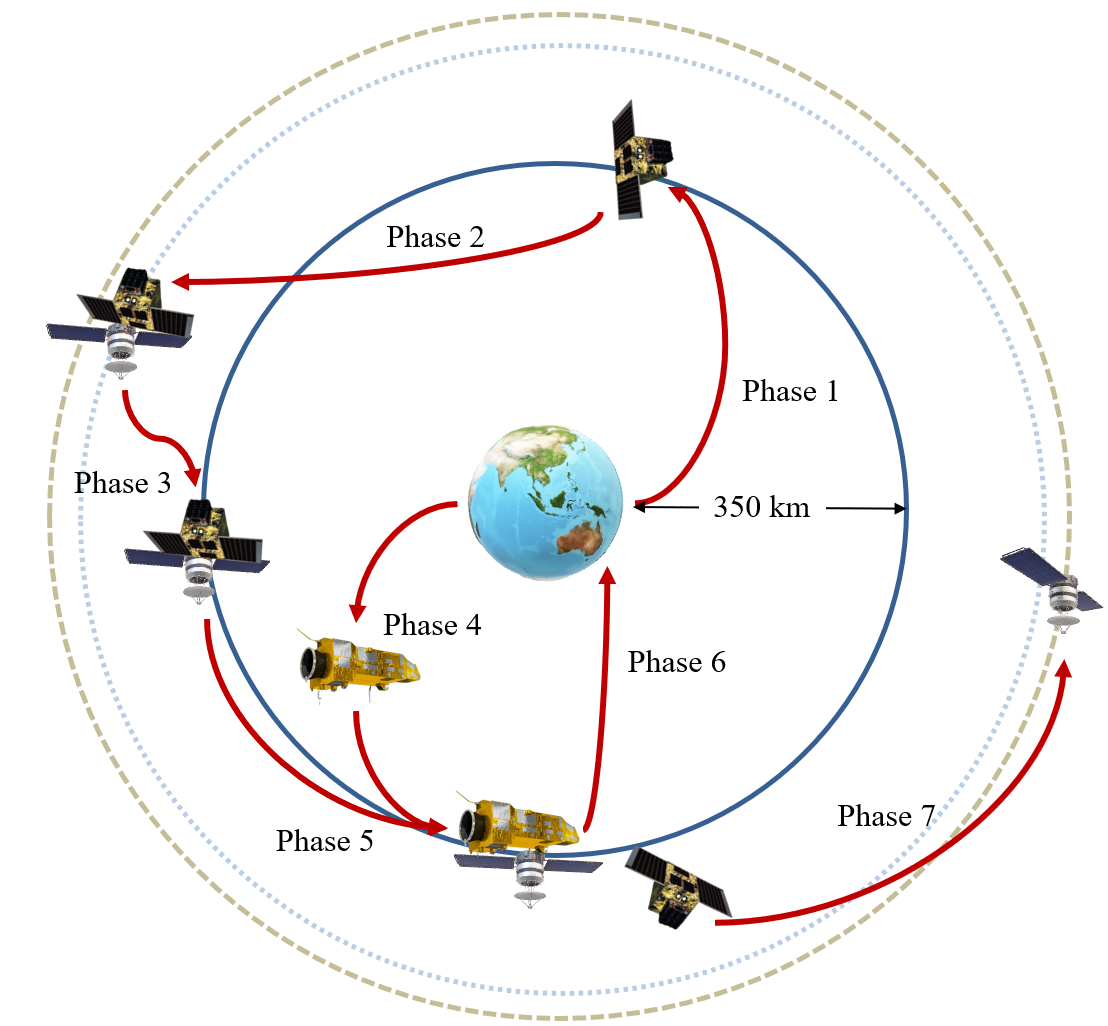}
\end{minipage}
\begin{minipage}[b]{0.3\linewidth}
\begin{itemize}
    \item Phase 1: Servicer launched 
    \item Phase 2: Servicer docks and detumbles debris 
    \item Phase 3: Servicer brings debris to a low altitude
    \item Phase 4: Reentry Shepherd launched
    \item Phase 5: Servicer holds on to debris while Shepherd docks
    \item Phase 6: Shepherd performs controlled reentry with debris 
    \item Phase 7: Servicer returns to collect the next debris
\end{itemize}
\end{minipage}
    \caption{Mission architecture of the multi-ADR mission}
    \label{ADR}
\end{figure}

The proposed mission architecture of the multi-ADR mission is the same as discussed in \cite{ADRMW}, and is shown in Figure \ref{ADR}. In this mission, two spacecraft are involved in the debris removal. The targeted debris are tonne class objects, including discarded rocket bodies and large derilict satellites \footnote{This mission was proposed as part of a collaborative study between the University of Auckland, Astroscale and Rocketlab.}.\par 

First, a Servicer spacecraft approaches and rendezvous with the target debris. When the rendezvous is achieved, the Servicer brings the object to a low-altitude orbit at $350$~km. The debris is then handed over to another spacecraft - named the Reentry Shepherd- which docks with the debris and performs controlled reentry on its behalf. A handover altitude of 350~km was selected to reduce the $\Delta v$ required by the Servicer by minimizing the orbital transfers it needs to perform while ensuring the technical feasibility and satisfying safety constraints posed by the altitude of the International Space Station \cite{astroscale}. \par

The Servicer can be reused for several debris removals, while each Reentry Shepherd can only be used once as it burns while deorbiting the debris. The proposed mission can perform multi-ADR services cheaper than a single spacecraft, which would perform all the mission phases and then burn in the atmosphere after removing a~single tonne-class debris \cite{astroscale}. This is because a space mission's development and operation costs are proportional to the system's dry mass~\cite{Jones2015}. Thus, although launching one Servicer and $n$ Shepherds to remove $n$ debris requires one more launch than using $n$ individual spacecraft, it is expected to be cheaper because of the lower overall mass launched. \par 

In this paper, guidance for the servicer spacecraft going from the initial 350 km orbit to the debris- called up leg- and coming back to the 350 km with the debris- called down leg- is studied. 

\section{MPC Guidance Methodology}
\begin{figure}[hbt!]
    \centering
    \includegraphics[width = 0.9\textwidth]{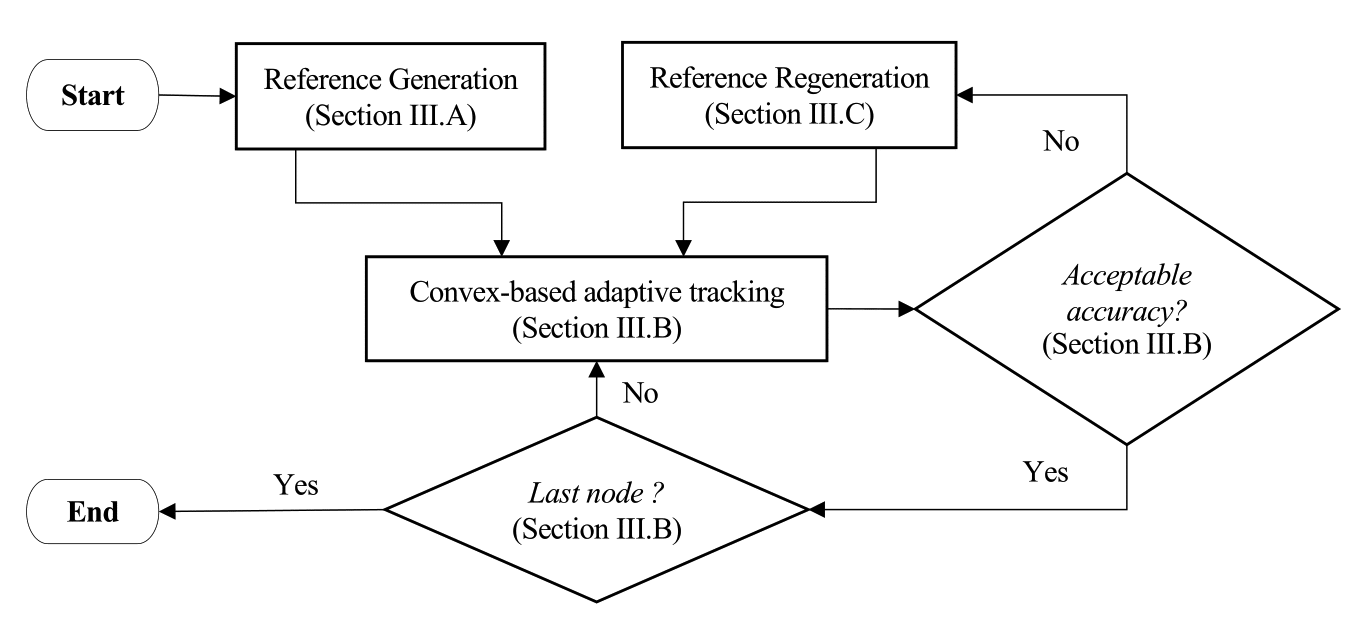}
    \caption{High-level overview of the MPC guidance}
    \label{1}
\end{figure}

The proposed guidance methodology is summarized in Figure \ref{1}, where the MPC architecture is divided into three stages.  
\begin{enumerate}
    \item \textbf{Reference Generation (Section \ref{Refgen})} A reference trajectory and a control profile to go from a mean initial semi-major axis $\bar{a}_0$, inclination $\bar{i}_0$, and RAAN $\bar{\Omega}_0$ to a mean target semi-major axis $\bar{a}_T$, inclination $\bar{i}_T$, and RAAN $\bar{\Omega}_T$ are determined using the PMDT. To incorporate phase angle matching into the transfer, a phase angle profile that reaches the targeted phase is also computed.
    \item \textbf{Convex-based Adaptive Tracking (Section \ref{convextracking})} The PMDT reference is segmented into equal time segments and tracked using a single iteration convex-optimization scheme. Then, the optimized control is forward propagated under realistic, nonlinear dynamics that include thrust errors to simulate real spacecraft motion.  Following the forward propagation, if phase matching is desired, the endpoint of each tracking segment is shifted in time, such that the real trajectory can follow the desired phase angle profile from Section \ref{Refgen}. 
    \item \textbf{Reference Regeneration (Section \ref{RefRegen})}: The PMDT reference is recomputed if the spacecraft deviates from the original reference significantly (i.e, if the accuracy is unacceptable). 
\end{enumerate}

 \tikzstyle{connector} = [draw, -latex']
 \tikzstyle{terminator} = [rectangle, draw, text centered, rounded corners, minimum height=2em]
 \tikzstyle{process} = [rectangle, draw, text centered, minimum height=2em]
 \tikzstyle{decision} = [diamond, draw, text centered, minimum height=2em]

\begin{figure}[hbt!]
    \centering
    \includegraphics[width = \textwidth]{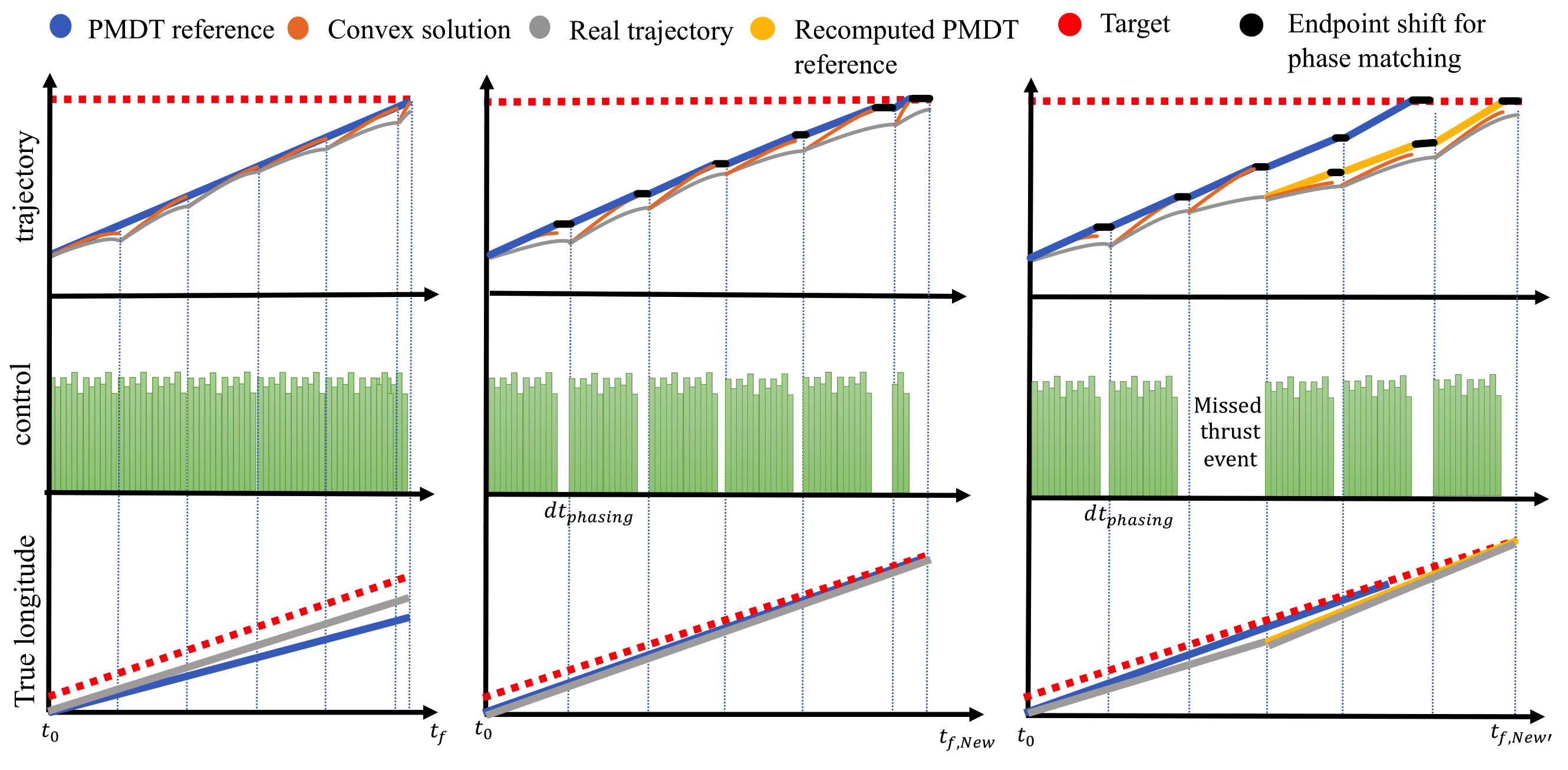}
    \caption{Overview of the convex-based  {MPC} (left: without reference recomputations and phase matching, center: with phase matching, right: with phase matching and reference recomputation following a misthrust event)}
    \label{MPCprocess}
\end{figure}

Figure \ref{MPCprocess} provides a graphical illustration of these steps. Figure \ref{MPCprocess} (left) shows the case where no phase matching is done and no recomputations are required to reach the target.  Figure \ref{MPCprocess} (middle) illustrates a case where phase matching is conducted but no recomputations are needed. Lastly, Figure \ref{MPCprocess} (right) shows a case with phase matching, where a recomputation is needed to reach the target due to a missed thrust event.    \par

\subsection{Reference Generation}\label{Refgen}

\begin{figure}[hbt!]
    \centering
    \includegraphics[width = 0.9\textwidth]{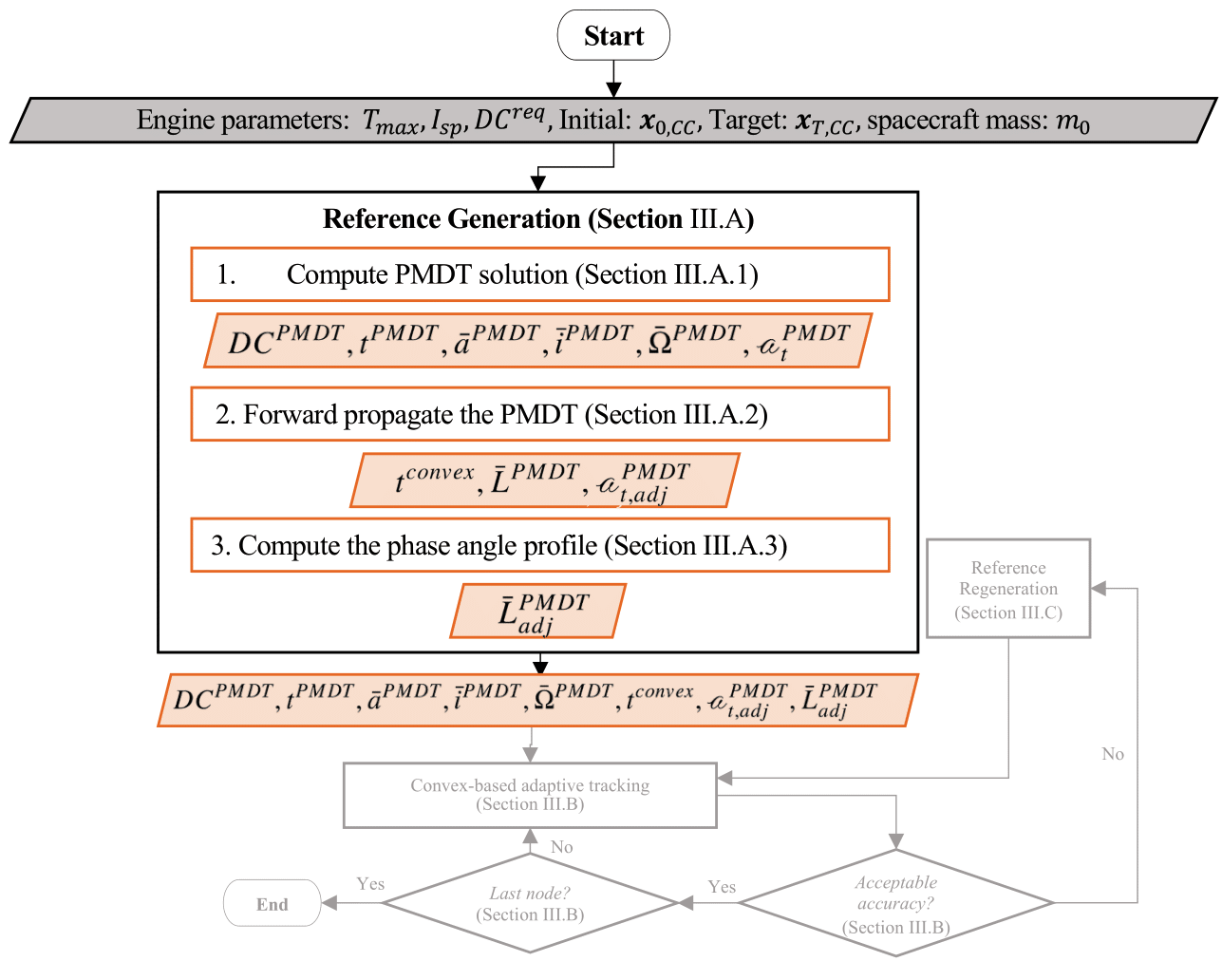}
    \caption{Reference generation process}
    \label{Refgenfig}
\end{figure}

The first step of the MPC guidance is the reference generation. Figure \ref{Refgenfig} illustrates the substeps involved. The inputs to the reference generation include the maximum spacecraft thrust $T_{max}$, specific impulse $I_{sp}$, initial spacecraft wet mass $m_0$,  duty cycle $DC^{req}$, and launch date $t_0$, along with the mean initial coordinates $\bar{a}_0, \bar{i}_0, \bar{\Omega}_0$  and mean target coordinates $\bar{a}_T, \bar{i}_T, \bar{\Omega}_T$.   \par 
Note that the inputs and outputs of the reference generation are mean coordinates given in the true-of-date (TOD) frame. When necessary, the mean to osculating- and vice versa- conversions are conducted as shown in \cite{8}. In this paper, all mean coordinates will be denoted by an overbar ($\ \bar{} \ $), any state without an overbar must be considered as osculating. 
The conversion to and from the J2000 frame to TOD requires accounting for precession and nutation and is done as given in \cite{DEconv}.  
Hence, if the initial and target coordinates are given in osculating, cartesian elements as $\bm{x}_{0,CC}$ and $\bm{x}_{T,CC}$, they must be converted to mean Classical Orbital Elements (COE), $\bar{\bm{x}}_{0,COE} = [\bar{a}_0, \bar{e}_0, \bar{i}_0, \bar{\Omega}_0, \bar{\omega}_0, \bar{\theta}_0]$ and $\bar{\bm{x}}_{T,COE} = [\bar{a}_T, \bar{e}_T, \bar{i}_T, \bar{\Omega}_T, \bar{\omega}_T, \bar{\theta}_T]$ prior to running the reference generation.  The outputs of each step of the reference generation will be elaborated in the corresponding sections.

The reference generation is divided into three substeps.
\begin{enumerate}
    \item \textbf{Computation of the PMDT solution (Section \ref{PMDT})} The PMDT discussed in \cite{ADRMW} is used to generate a reference for the guidance.  A grid search is added to avoid search space discontinuities, and a thrust margin is implemented via the duty cycle such that the guidance shall have more thrust capability in reality than in the reference. 
    \item \textbf{PMDT forward propagation (Section \ref{fPMDT})} The PMDT solution is forward propagated, and another margin of thrust is added if the forward propagated undershoots the target in reality. 
    
    \item \textbf{Phase angle profile calculation (Section \ref{Ladjcal})}: A phase angle profile to be tracked is calculated based on the forward propagation from the previous section, such that phase matching can be conducted alongside the transfer when necessary. 

\end{enumerate}

\subsubsection{PMDT Computation}\label{PMDT}
The PMDT- first introduced in \cite{ADRMW}- can generate a fuel or time optimal trajectory to go from an initial, mean, TOD-frame $\bar{a}_0, \bar{i}_0, \bar{\Omega}_0$ to target $\bar{a}_T, \bar{i}_T, \bar{\Omega}_T$ through local optimization methods. This is done by selecting an optimal drift orbit- mean semi-major axis $\bar{a}_d$ and inclination $\bar{i}_d$- that utilizes the effect of $J_2$ to reach the $\bar{\Omega}_T$, to perform the RAAN change without any fuel consumption.  It consists of two main algorithms, the Extended Edelbaum method (Algorithms \ref{A1}) and the RAAN matching scheme (Algorithms \ref{A2}).

\begin{figure}[hbt!]
\centering
      \includegraphics[width = 0.7\textwidth]{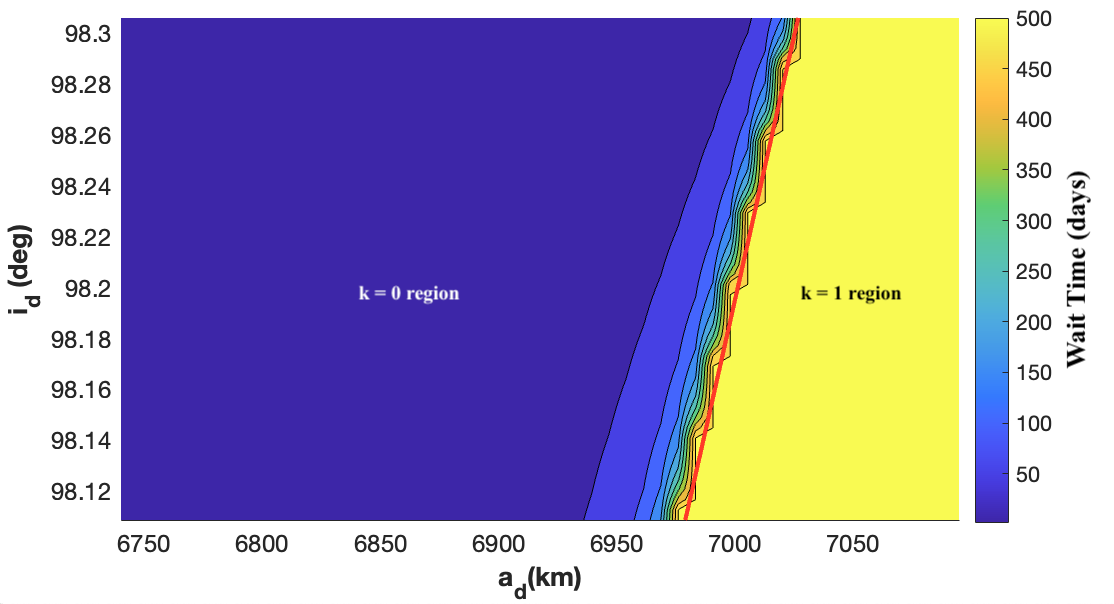}
\caption{The discontinuous search space encountered during RAAN matching (red line: Differential drift rate between the target and drift orbit equals zero) }
    \label{issue1}
  \label{fig:test1}
  \end{figure}

\begin{itemize}
\item \textbf{Extended Edelbaum method (Algorithm \ref{A1})} is a version of the classical Edelbaum method adapted to consider the effect of atmospheric drag and duty cycles. This method can generate a $TOF$ and $\Delta v$ optimal trajectory to reach a desired semi-major axis and inclination. To consider the effect of the $DC$, the maximum thrust $T_{max}$ is multiplied by it. $\bar{\Omega}$ is updated to consider the effect of $J_2$, however, the impact of the thrust on $\bar{\Omega}$ is neglected in this Algorithm.

\item \textbf{RAAN matching scheme (Algorithm \ref{A2})} builds on Algorithm \ref{A1} such that RAAN changing transfers can be optimized. In this method, $J_2$ is used to achieve a target RAAN by drifting at an intermediate orbit, as done in \cite{6}. This results in a transfer that has a thrust-coast-thrust structure, where the thrust arcs are computed using the above Extended Edelbaum method. The drift orbit variables $\bar{a}_d$ and $\bar{i}_d$ are obtained by optimizing the transfer for time or propellant consumption using a local, gradient-based optimization algorithm.  When drifting, thrust is utilized to counteract the effect of drag. Hence, the thrust during drift is set to be equal in magnitude and opposite in direction to the drag acceleration experienced.  
The outputs of the RAAN matching scheme are the optimal $TOF^{PMDT}$,  $\Delta v^{PMDT}$, thrust profile $\mathscr{a}^{PMDT}_t$ and the optimal trajectory profile [$\bar{a}^{PMDT}, \bar{i}^{PMDT}, \bar{\Omega}^{PMDT}$] defined in time vector $t^{PMDT}$. 
\end{itemize}
In the RAAN matching scheme,  the drift time is obtained from Eq.\eqref{tw}, where it is a function of the number of RAAN revolutions (+$2 k \pi$ where $k\in \mathbb{N}$ ). As such, there are discontinuous jumps in the wait time as the integer $k$ varies.  This effect is demonstrated in Figure \ref{issue1}, where a jump in drift time is seen between $k =0$ and $k = 1$.   It was noted that these discontinuities cause convergence problems for gradient-based local optimization methods such as the interior point method used in the RAAN matching scheme.  Hence, in this work, an initial guess for the drift orbit was first obtained through a grid search that finds a minimum $\Delta v/TOF$ solution that is away from any discontinuities. This was given as a guess for the interior-point method used in the RAAN matching scheme such that the optimization steers clear of the discontinuous regions in the search space. \par 

\begin{figure}[hbt!]
\centering
\begin{minipage}{.48\textwidth}
  \centering
      \includegraphics[width = \textwidth]{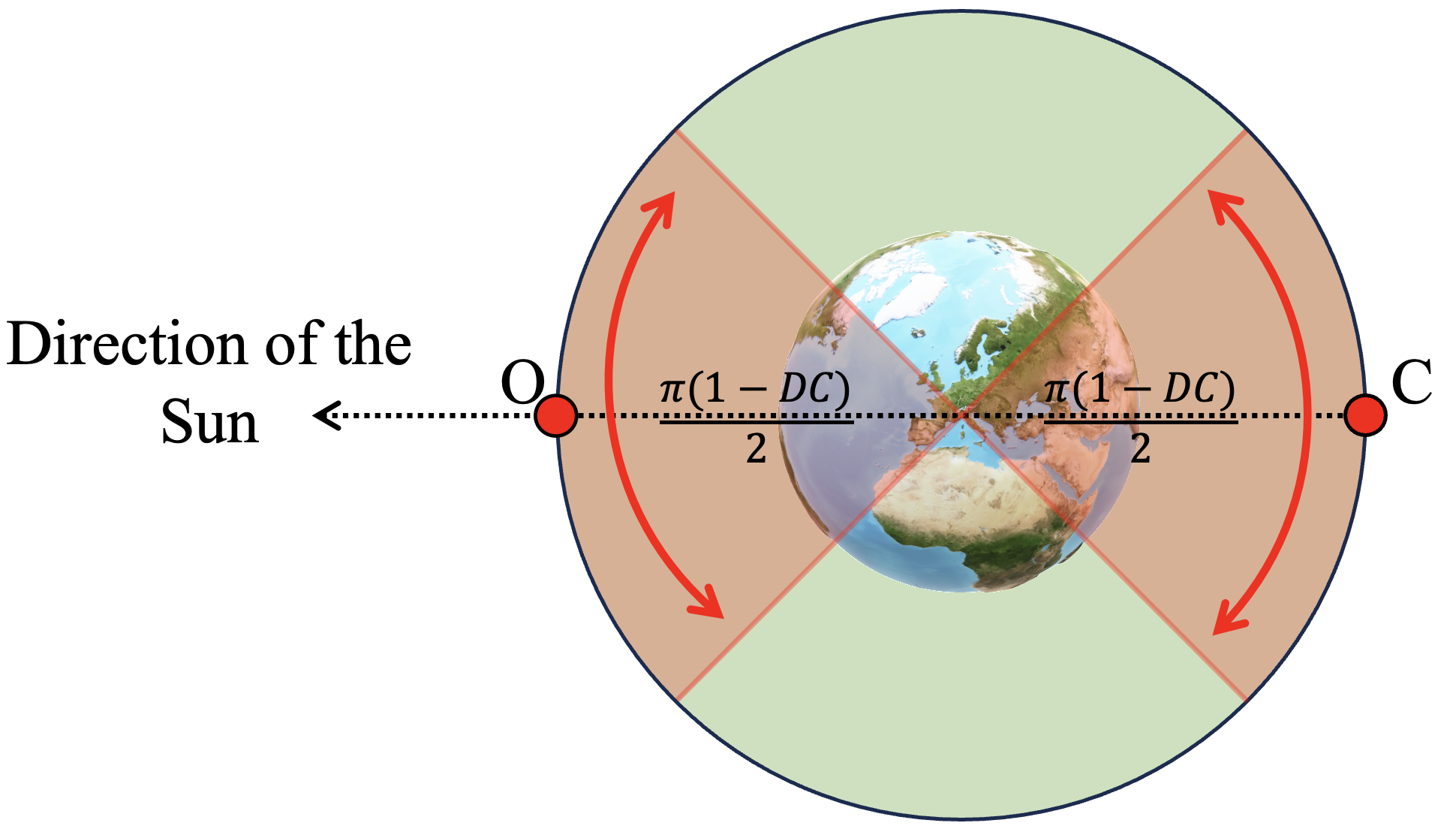}
\caption{The PMDT thrust profile formulation \cite{ADRMW}: the coast arcs are shown in red, while the thrust arcs are shown in green. }
    \label{eclprop}
\end{minipage}%
\hspace{0.01\textwidth} 
\begin{minipage}{0.48\textwidth}
  \centering
\includegraphics[width = \textwidth]{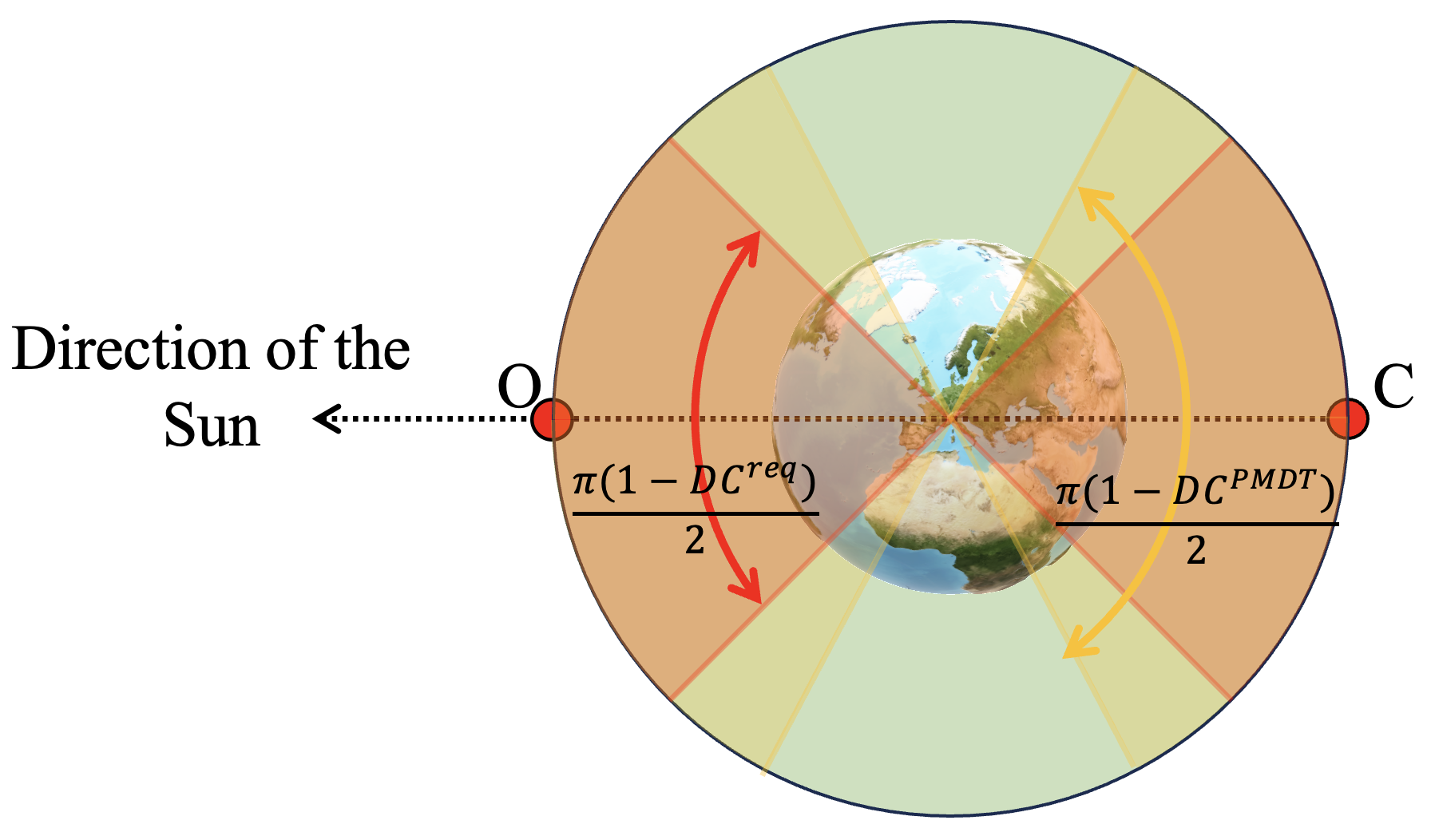}
    \caption{The difference in thrust profile for $DC^{PMDT}$ (reference duty cycle (yellow coast arc)) and $DC^{req}$ (required duty cycle (red coast arc)).}
    \label{issue2}
\end{minipage}
\end{figure}

In this work and in \cite{ADRMW}, the effect of eclipses and $DC$ are combined in the PMDT by turning the thrust off symmetrically across the orbit, centered around the center of the eclipse ($C$) and its antipodal counterpart ($O$), such that the total thrust duration agrees with the duty cycle. This is shown in Figure \ref{eclprop}. This symmetrical coasting is conducted to minimize the eccentricity build-up \cite{Viavattene2022DesignOM}. When generating a reference for the guidance, the PMDT is calculated at a $DC$ value lower than the required one. The duty cycle for the reference calculation is denoted $DC^{PMDT}$, while the required spacecraft duty cycle is denoted $DC^{req}$. Having $DC^{PMDT} < DC^{req}$ implies that the coast arcs of the spacecraft are shorter in reality than in the reference,  as shown in Figure \ref{issue2}. This retains a margin of thrust that could be used to handle any mismatch between the designed trajectory and the guidance that arises from thrust uncertainties and orbital perturbations that are not accounted for by the PMDT reference.

The following outputs from the PMDT are necessary for other steps of the convex-MPC guidance: the time vector $t^{PMDT}$, mean semi-major axes profile $\bar{a}^{PMDT}$, mean inclination profile $\bar{i}^{PMDT}$, mean RAAN profile $\bar{\Omega}^{PMDT}$ and thrust acceleration profile $\mathscr{a}_t^{PMDT}$.

\subsubsection{Forward Propagation of the PMDT}\label{fPMDT}
The PMDT trajectory is forward propagated to determine if it is able to reach the target.
The forward propagation requires $DC^{PMDT}, t^{PMDT}$ and $\mathscr{a}_t^{PMDT}$ parameters from the previous section, along with the spacecraft parameters $T_{max}$, $I_{sp}$, $m_0$, and $\bm{x}_{0,CC}$ (initial osculating cartesian state of the spacecraft in the TOD frame). The time vector for the forward propagation is set such that
$t^{convex}  = t_0 : P/N: t_f$,
where $P, N$ and $t_f$ represent the orbital period at the start, the number of nodes per orbit set for the convex optimization, and the end time of the mission, respectively.  The propagation algorithm is given as Algorithm 3 in \cite{ADRMW} (provided as Algorithm \ref{alg:2} in the Appendix).  The main outputs of Algorithm \ref{alg:2} are the final state reached $\bm{x}_{f,CC}$, given in osculating cartesian coordinates in TOD frame, and the profile of the mean argument of latitude $\bar{L}^{PMDT}$.\par
Following the forward propagation, $\bm{x}_{f,CC}$ is converted to mean COEs in TOD frame to provide  $\bar{\bm{x}}_{f,COE} = [\bar{a}_f, \bar{e}_f, \bar{i}_f, \bar{\Omega}_f, \bar{\omega}_f, \bar{\theta}_f]$. Then, if $\bar{a}_f, \bar{i}_f,\bar{\Omega}_f$ undershoots the targeted $\bar{a}_T, \bar{i}_T, \bar{\Omega}_T$, an additional margin of thrust is added to aid the guidance. This addition is done by first estimating the $\Delta v$ to go from $\bar{a}_f, \bar{i}_f,\bar{\Omega}_f$ to $\bar{a}_T, \bar{i}_T,\bar{\Omega}_T$ using Eq. \eqref{dvp1}. Note that this equation is obtained by integrating the Gauss Variational Equations of maximum rates of change given in \cite{GVE} in time.

\begin{align}
     \bm{\Delta v}'&= [{\Delta v_{\bar{a}}}, {\Delta v_{\bar{i}}},{\Delta v_{\bar{\Omega}}}]^T  = \begin{bmatrix}  \frac{\Delta \bar{a}}{2 }\sqrt{\frac{\mu (1-\bar{e}_f)}{\bar{a}_f^3(1+\bar{e}_f)}} \\ \sqrt{\frac{\mu}{\bar{a}_f(1 - \bar{e}_f^2)}} \Delta \bar{i} (\sqrt{1 - \bar{e}_f^2 \sin^2 \bar{\omega}_f} -\bar{e}_f |\cos \bar{\omega_f} | \\   \sqrt{\frac{\mu}{\bar{a}_f(1 - \bar{e}_f^2)}} \Delta \bar{\Omega} \sin{\bar{i}} (\sqrt{1 - \bar{e}_f^2 \cos^2 \bar{\omega}_f} - \bar{e}_f |\sin \bar{\omega}_f |) 
  \label{dvp1}
  \end{bmatrix}
\end{align}
 
where $\Delta \bar{a} = \bar{a}_f - \bar{a}_T$ , $\Delta \bar{i} = \bar{i}_f - \bar{i}_T$, and $\Delta \bar{\Omega} = \bar{\Omega}_f - \bar{\Omega}_i$. The PMDT $\Delta v$ profile ${\Delta v^{PMDT}}$ is then adjusted to include $\lVert \bm{\Delta v}' \rVert $ as

\begin{equation}
    \Delta v^{PMDT}_{adj}(t^{PMDT}) = \Delta v^{PMDT}(t^{PMDT}) + \frac{t^{PMDT}}{t_f-t_0}\lVert \bm{\Delta v}' \rVert 
\end{equation}
to ensure a uniform increment of $\Delta v$. Then the new spacecraft mass profile is calculated using the classical rocket equation
     
\begin{equation}
m^{PMDT}_{adj}(t^{PMDT}) = {m_0}/{\exp \left( \frac{ \Delta v^{PMDT}_{adj}(t^{PMDT}) }{I_{sp}g_0} \right)}, 
\end{equation}
     
following which an augmented thrust profile can be obtained as 
 
\begin{equation}
  \mathscr{a}^{PMDT}_{t,adj}(t) = \frac{T_{max}}{m^{PMDT}_{adj}(t)}.
\end{equation}

The following outputs from the PMDT forward propagation are necessary for other steps of the convex-MPC guidance: the adjusted acceleration profile $\mathscr{a}^{PMDT}_{t,adj}(t)$, time vector $t^{convex}$ and the true longitude profile $\bar{L} ^{PMDT}$.

\subsubsection{Computing the Phase Angle Reference}\label{Ladjcal}
 
For the up leg of the mission, the spacecraft is expected to rendezvous with the target debris. As such, the phase of the spacecraft must be matched with that of the target, along with other orbital elements. In this work, the phasing angle is the true longitude. The forward propagation in  Algorithm \ref{alg:2} provides an estimate of the mean true longitude profile $\bar{L} ^{PMDT}$ of the spacecraft. Here, the slope of the $\bar{L} ^{PMDT}$ is modified as follows such that it reaches the mean true longitude of the debris at $t_f$, denoted as $\bar{L}_T (t_f)$.

\begin{equation}
    \bar{L} ^{PMDT}_{adj}(t^{convex}) = \bar{L} ^{PMDT}(t^{convex}) +   \frac{t^{convex}}{t_f-t_0} {\delta \bar{L}}
\end{equation}

where $ \delta \bar{L} = \bar{L}_T (t_f) - \bar{L} ^{PMDT}(t_f)$. The calculated $\bar{L} ^{PMDT}_{adj}$ is then given as input to the convex-based tracking to follow. \par 
Note that when selecting a duty cycle $DC^{PMDT}$ for the reference, it must be ensured that $\delta \bar{L}$ is close to zero. Otherwise, the convex tracking would experience large shifts in the segment end times, making the reference invalid soon after the departure. Through experimentation, it was realized that $\delta \bar{L}$ must be at least less than 45 deg to ensure adequate tracking. 
Hence, it is ideal to run the PMDT reference and forward propagated it for a range of $DC^{PMDT}$ values that are smaller than $DC^{req}$, and then select the $DC^{PMDT}$ for which the $\delta \bar{L}$ is the smallest as the reference $DC^{PMDT}$. \par 
The following output from the phase angle reference calculation is necessary for other steps of the convex-MPC guidance: the adjusted phase angle profile $ \bar{L} ^{PMDT}_{adj}$. \par 

\subsubsection*{Outputs of Reference Generation}
The following outputs from the reference generation process are required  
to proceed to convex tracking:  the mean semi-major axes profile $\bar{a}^{PMDT}$, mean inclination profile $\bar{i}^{PMDT}$, mean RAAN profile $\bar{\Omega}^{PMDT}$, and the time vector $t^{PMDT}$ from Section \ref{PMDT},  the adjusted thrust acceleration profile $\mathscr{a}^{PMDT}_{t, adj}$ and the time vector of node placement $t^{convex}$ from Section \ref{fPMDT}, and the adjusted thrust profile $\bar{L}^{PMDT}_{adj}$ from Section \ref{Ladjcal}.

\subsection{Convex-based Adaptive Tracking} \label{convextracking}

\begin{figure}[hbt!]
    \centering
    \includegraphics[width = 0.9\textwidth]{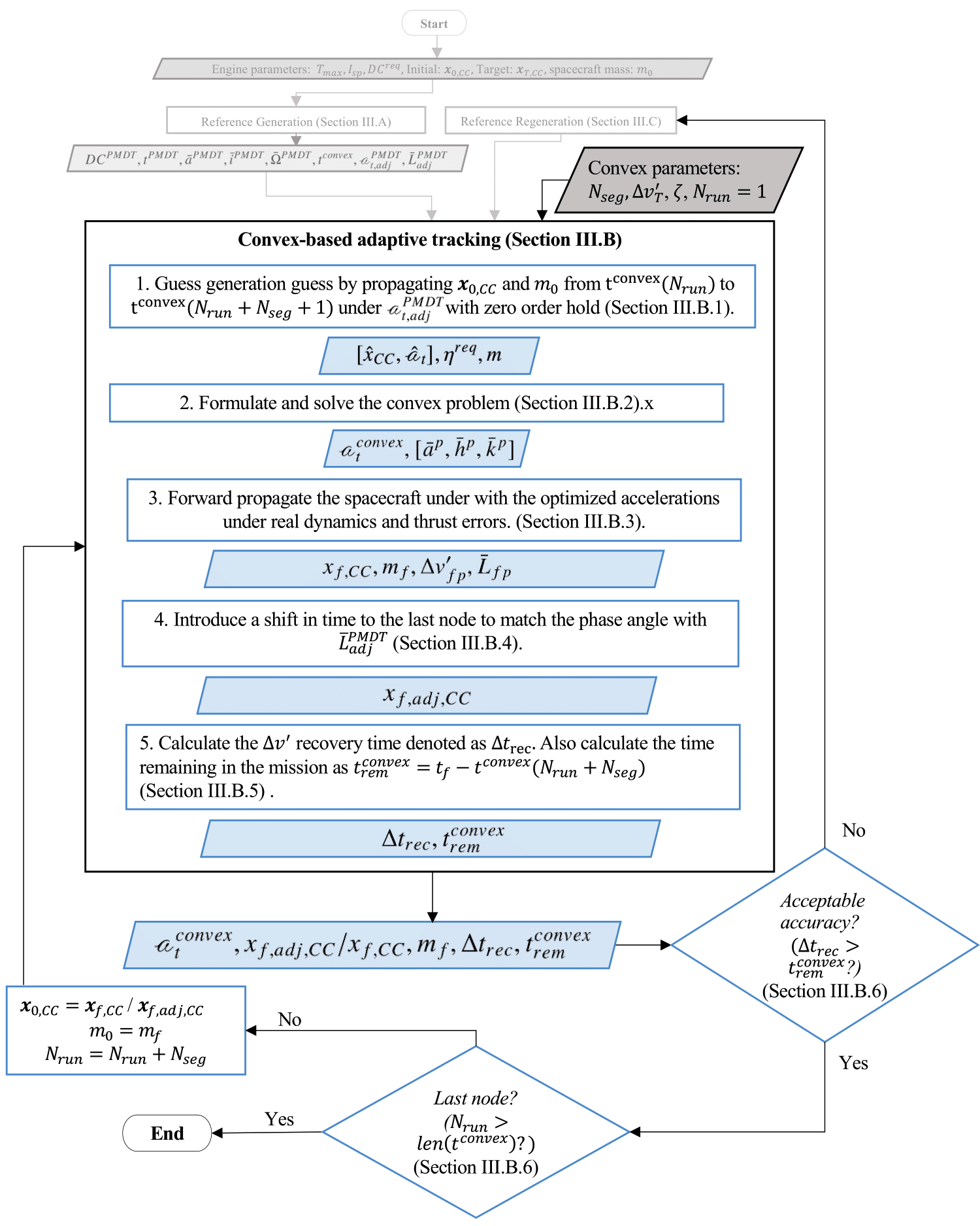}
    \caption{Convex-based adaptive tracking algorithm}
    \label{2}
\end{figure}

Following the reference generation, convex-optimization-based adaptive tracking is used to guide the spacecraft along the reference. This is an iterative process that involves segment-wise initial guess generation (Section \ref{initialg}), convex optimization (Section \ref{convex}), and forward propagation (Section \ref{fprops}) till the end of the transfer (i.e., the end of the time vector $t^{convex}$) is reached. The high-level overview of this process is given in Figure \ref{2}, and the outputs of each step will be elaborated in the corresponding sections. Note that the convex tracking is also performed in the TOD reference frame. 

Firstly, the time vector $t^{convex}$ is split into $N_{seg}$ long segments where
$N_{seg} = nN$. $n$ is the number of orbits per tracking segment and $N$ is the number of convex nodes per orbit. Both $n$ and $N$ are user-set parameters. At the start of the convex tracking, the starting osculating state is denoted $\boldsymbol{x}_{0,CC}$, and the starting mass is $m_0$.  The index $N_{run}$ starts at one and indexes through the vector $t^{convex}$. $\zeta$ is another user-set parameter related to the convex optimization and will be discussed in Section \ref{convex}.  
\par 
Following an iteration of the convex tracking (i.e., steps 1-5 in Figure \ref{2}), an accuracy check is conducted to see if the endpoint of the forward propagated trajectory is sufficiently close to the PMDT reference at that time (Section \ref{Acheck}). If so, a second check is conducted to see if the iteration has reached the last node (Section \ref{Lcheck}). If the tracking is not at the last node, it continues on to the next segment.  If the accuracy check is deemed unacceptable, a new reference is calculated (Section \ref{RefRegen}).

\subsubsection{Initial Guess Generation} \label{initialg}
An initial guess is generated for a single segment from $t^{convex}(N_{run})$ to $t^{convex}(N_{run} + N_{seg}+1)$ by forward propagating the adjusted thrust accelerations $\mathscr{a}_{t,adj}^{PMDT}$-starting from $\bm{x}_{0,CC}$, as given in Algorithm \ref{alg:3}. Inputs for Algorithm \ref{alg:3} are the initial osculating Cartesian coordinates $\bm{x}_{0,CC}$, starting mass $m_0$, convex time vector $t^{convex}$,  PMDT time vector $t^{PMDT}$ and the PMDT adjusted acceleration profile $\mathscr{a}_{t,adj}^{PMDT}$.

\begin{algorithm}[hbt!]
\caption{Initial guess generation for a convex tracking segment}\label{alg:3}
\textbf{Input}  $\bm{x}_{s, CC} = [\bm{r}_s, \bm{v}_s]^T = \bm{x}_{0,CC}$, $m_0$, $t^{convex}$, $t^{PMDT}$, $  \mathscr{a}^{PMDT}_{t,adj}$
\setstretch{0.75}
\begin{algorithmic}
\State Define the time vector for the segment. $t = t^{convex}(N_{run} : N_{run} + N_{seg}) $ and let $m(1) = m_0$. 
\State Calculate the mean COE coordinates.  $\bar{\bm{x}}_{s,COE} = oscCart2meanKep(\boldsymbol{x}_{s,CC})$. \Comment{conversion given in \cite{8}}
\For{$i = 1: N_{seg}$}
\State  Calculate the mean argument of latitude: $\bar{L}(t) =  \bar{\bm{x}}_{s,COE}(5) + \bar{\bm{x}}_{s,COE}(6) $
\State  Calculate the eclipse/DC profiles corresponding to $DC^{req}$ and $DC^{PMDT}$ as given in \cite{ADRMW}.
\State \hspace{\algorithmicindent}    Calculate $q_1 = \cos^{-1}(\cos(\bar{L} - \bar{L}_c)) $  and $q_2 =  \cos^{-1}(\cos(\bar{L}- \bar{L}_c -\pi )) $ \par \Comment{$\bar{L}_c$ is the center of the eclipse, calculated as discussed in \cite{ADRMW}}
 \State \hspace{\algorithmicindent}    \textbf{if} { $q_1 < \frac{\pi}{2}(1-DC^{PMDT})$ \textbf{or}  $q_2 < \frac{\pi}{2}(1-DC^{PMDT})$}  \textbf{then} $\eta^{PMDT}(i)=0$ 
  \textbf{else} $\eta^{PMDT}(i)=1$ 
 \State \hspace{\algorithmicindent}    \textbf{if} { $q_1 < \frac{\pi}{2}(1-DC^{req})$ \textbf{or}  $q_2 < \frac{\pi}{2}(1-DC^{req})$}  \textbf{then} $\eta^{req}(i)=0$ 
  \textbf{else} $\eta^{req}(i)=1$  \Comment{$\eta(i)=0$ indicates drifting.}
\State Calculate a guess for the control acceleration magnitude $ \hat{\mathscr{a}}^{PMDT}_t(i)$ by interpolating and averaging  $\mathscr{a}_t^{PMDT}$.
\begin{align}
 \hat{\mathscr{a}}_t(i) &= \frac{\eta^{PMDT}(i) }{2DC^{PMDT}} [\textrm{interp}(t^{PMDT} , \mathscr{a}^{PMDT}_{t,adj} , t(i)) +  (\textrm{interp}(t^{PMDT} , \mathscr{a}^{PMDT}_{t,adj} , t(i+1)) ] \end{align}
 \State Calculate $\Delta v(i) =   \hat{\mathscr{a}}_t(i) (t(i+1) -t(i))$.
 \State Calculate mass after the burn: $m(i+1) = {m(i)}/{\exp(\Delta v/I_{sp}/g_0)}$ \Comment{The mass must be kept track of as the drag acceleration $\bm{\mathscr{a}}_{drag}$ is a function if it.}
\State Determine the average out-of-plane angle ($\beta$) for the burn. 
\State \hspace{\algorithmicindent}  Calculate the $\beta$ angles at $t(i)$ and $t(i+1)$  using Eq. 6 in \cite{ADRMW} and average. $\beta_{avg} = \frac{1}{2}(\beta(i)  + \beta(i+1))$
\State  Reverse thrust direction across the line of nodes to make the inclination change possible.
 \If{$\bar{i}_f > \bar{i}_0$}
    \State  \textbf{if}  { $\bar{L}(t)   > \pi/2$ or $\bar{L}(t)  \leq 3\pi/2$} \textbf{then} $\beta = -|\beta| $ \textbf{else}  { $\beta = |\beta|$}
  \Else
    \State  \textbf{if} { $\bar{L}(t)  \geq \pi/2$ or $\bar{L} (t)  < 3\pi/2$} \textbf{then}
$\beta = |\beta|$  \textbf{else} 
$\beta = -|\beta|$
  \EndIf
\State Calculate the RTN acceleration:  { $ \hat{\bm{\mathscr{{a}}}}_{t,RTN} = \hat{\mathscr{a}}_t(i) [0 , \cos{\beta_{avg}}, \sin{\beta_{avg}}]$ and store $\hat{\bm{x}}_{CC}(i,:) = \boldsymbol{x}_s, \hat{\bm{\mathscr{a}}}_{t}(i,:) = \hat{\bm{\mathscr{a}}}_{t,RTN}$.}
\State  {Propagate from $t(i \rightarrow i+1)$ with $\bm{\mathscr{a}}_{g}, \bm{\mathscr{a}}_{J_2}. \bm{\mathscr{a}}_{drag}$ (calculated as in \cite{8}) and ${\hat{\bm{\mathscr{a}}}_{t,TOD}}$ accelerations.}  \Comment{$ \hat{\bm{\mathscr{a}}}_{t,RTN}$ must be converted to TOD for propagation}
 
\State  { Obtain the state at $t(i+1)$ and define it as the new $\bar{\bm{x}}_{s,CC}$.} 
\State  Obtain the new mean COE coordinates. $\bar{\bm{x}}_{s,COE} = oscCart2meanKep(\boldsymbol{x}_{s,CC})$.
\EndFor
 \State Let $\hat{\bm{x}}_{CC} (N_{seg}+1,:) = \boldsymbol{x}_{s,CC}$, $\hat{\bm{\mathscr{a}}}_{t}(N_{seg}+1,:) = [0,0,0]^T$ \Comment{No thrust is allocated for the last node.}
\end{algorithmic}
\textbf{Output} $\hat{\bm{x}}_{CC}$, $\hat{\bm{\mathscr{a}}}_{t}$, $\eta^{req}, m $.
\end{algorithm}

 Note that Algorithm \ref{alg:3} implements a zero-order hold on the RTN acceleration obtained from the PMDT, such that the initial guess closely resembles the convex formulation, as it also utilizes a zero-order hold on control at each node.  Also, note that a segment of convex tracking has $N_{seg}+1$ nodes, but the last node has no associated control.  \par 
 The following outputs from the initial guess generation are necessary for other steps of the convex-MPC guidance: the initial guess of the state and control $[\hat{\bm{x}}_{CC}, \hat{\bm{\mathscr{a}}}_{t}]$, the thrust profile $\eta^{req}$ for the required spacecraft duty cycle $DC^{req}$, and the mass profile $m$.

\subsubsection{Formulating and Solving the Convex Optimization Problem}  \label{convex}
Once an initial guess is calculated for a segment, convex optimization is used to obtain the optimal acceleration profile that minimizes the total fuel consumption and the distance between the real trajectory and the PMDT reference at the end of the segment.  Numerically, this is denoted as 
\begin{equation}
    \textrm{minimize } J = \int^{t^{convex}(N_{run}+ N_{seg} )}_{t^{convex}(N_{run})} 	\lVert \bm{\mathscr{a}}_{t}^{convex}(t)\rVert dt  + \lVert \bm{\Delta v}' \rVert 
    \label{obj}
\end{equation}

where $\bm{\mathscr{a}}_{t}^{convex}(t) = [\mathscr{a}_r (t), \mathscr{a}_\theta (t), \mathscr{a}_n (t)]^T$ is the thrust acceleration profile in RTN coordinates. $\lVert \bm{\Delta v}' \rVert $ estimates the $\Delta v$ required to go from the mean, TOD orbital elements reached by the convex optimization given in modified equinoctial elements (MEEs) $\bar{\bm{x}}_{MEE}^{c} = [\bar{a}^{c}, \bar{f}^{c}, \bar{g}^{c}, \bar{h}^{c},\bar{k}^{c}, \bar{L}^{c}]$ to the PMDT reference orbital elements $[\bar{a}^{p}, \bar{h}^{p},\bar{k}^{p}] $ at $t^{convex}(N_{run} + N_{seg}+1)$, the endpoint of the segment being optimized.  The propagation of the dynamics is done in osculating elements; hence, an osculating to mean conversion using the formulae in \cite{8} is necessary to obtain $\bar{\bm{x}}^{c}_{MEE} $. As only $\bar{a},\bar{i}$ and $\bar{\Omega}$ are tracked, only the $\Delta v'$ due to $\bar{a},\bar{h}$ and $\bar{k}$ are included in $\lVert \bm{\Delta v}' \rVert $.  Firstly, $[\bar{a}^{p}, \bar{h}^{p},\bar{k}^{p}]$ is obtained by interpolating the PMDT output and converting them to MEEs. 
\begin{equation}
    x = \text{interp}(t^{PMDT}, x^{PMDT},T(N_{run} + N_{seg}+1)) \ \text{where} \ x = [\bar{a}, \bar{i}, \bar{\Omega}].
\end{equation}
Then, $ [\bar{a}^{p}, \bar{h}^{p},\bar{k}^{p}] = [\bar{a},\tan(\bar{i}/2)\cos{\bar{\Omega}}, \tan(\bar{i}/2)\sin{\bar{\Omega}} ]$. $\lVert \bm{\Delta v}' \rVert $ is obtained by integrating the modified equinoctial Gauss Variational Equations of maximum rates of change \cite{GVE} in time, which provides
    \begin{align}
     \bm{\Delta v'} = [{\Delta v_{\bar{a}}}, {\Delta v_{\bar{h}}},{\Delta v_{\bar{k}}}]^T = \left[ \frac{\Delta \bar{a}}{2 \bar{a}^{c} }\sqrt{\frac{\mu}{\bar{a}^{c}}}\sqrt{\frac{1-\bar{e}^p}{1+\bar{e}^p}},{2\Delta \bar{h}}\sqrt{\frac{\mu}{\bar{p}^{c}}} \frac{\sqrt{1-{\bar{g}^{c2}}}+\bar{f}^{c}}{{\bar{s}^{c2}}} ,{2\Delta \bar{k}}\sqrt{\frac{\mu}{\bar{p}^{c}}} \frac{\sqrt{1-{\bar{f}^{c2}}}+\bar{f}^{c}}{{\bar{g}^{c2}}} \right]^T
     \label{dvp2}
 \end{align}
 
 where $[\Delta \bar{a}, \Delta \bar{h}, \Delta \bar{k}] = [\bar{a}^c - \bar{a}^p, \bar{h}^c - \bar{h}^p, \bar{k}^c - \bar{k}^p]$. $\bar{e}^c= \sqrt{{\bar{f}^{c2}} + {\bar{g}^{c2}}}$ , $\bar{p}^c = \bar{a}^c(1- {\bar{e}^{c2}})$ and ${\bar{s}^{c2}} = 1 + {\bar{h}^{c2}} + {\bar{k}^{c2}}$.  Note that this is the same as $\Delta v'$ in Eq. \eqref{dvp1}, defined for MEEs. MEEs are used to implement this constraint to avoid the issue of equating angles a convex setting when matching the RAAN values. \par 
 
Note that the convex tracking is only required to match the reference $a, i$ and $\Omega$ precisely once the spacecraft reaches its target. Precise matching to the PMDT trajectory before reaching the target results in unnecessary fuel consumption. Hence, instead of using a hard boundary condition at the end of a tracking segment, the $\Delta v'$ from Eq. \eqref{dvp2} is minimized in the objective function. $\Delta v'$ also provides a singular accuracy measure, which is used to determine whether a recomputation of the reference is necessary in Section \ref{recT}. \par
 
 The convex tracking is also subjected to the following constraints and dynamics 

\begin{equation}
    \begin{aligned}
& \dot{\bm{x}}=\boldsymbol{f}(\boldsymbol{x}, \boldsymbol{u})= \bm{h(x)} + g(\bm{\mathscr{a}}_{t}^{convex}) \\
& \bm{x} \left(T(N_{run})\right)= \bm{x}_{0} \\
& \| \bm{\mathscr{a}}_{t}^{convex}\| \leq  \eta^{req} T_{max}/m \\
\end{aligned}
\label{cons}
\end{equation}

$\bm{h(x)}$ represents the natural dynamics, defined as $\bm{h(x)} = \bm{\mathscr{a}}_g +  \bm{\mathscr{a}}_{drag} + \bm{\mathscr{a}}_{J2}$ as given in Algorithm 1. $g: RTN \rightarrow TOD$ is a coordinate conversion function, which is implemented on $\bm{\mathscr{a}}_{t}^{convex}$  to convert it from RTN to TOD frame. $\eta^{req}$ was calculated in Section \ref{initialg} and represents the thrust profile for the segment based on the required duty cycle $DC^{req}$. $m$ is the mass profile for the segment, calculated in Section \ref{initialg}. 
\par 
{
The nonconvex nature of the objective function and problem dynamics in Eqs \eqref{obj} and \eqref{cons} implies that successive iterations may be required if the convexification error is high. However, using successive iterations is not ideal for real-time guidance, as it is time-consuming, and convergence to a global optimal is not guaranteed. As discussed in the Introduction, using the PMDT for reference generation increases the accuracy of the convexification. To retain and strengthen this convexification accuracy, a coordinate system that is considerably linear must chosen to represent the convex nodes.} \par
\begin{figure}[hbt!]
\centering
  \centering
    \includegraphics[width = 0.5\textwidth]{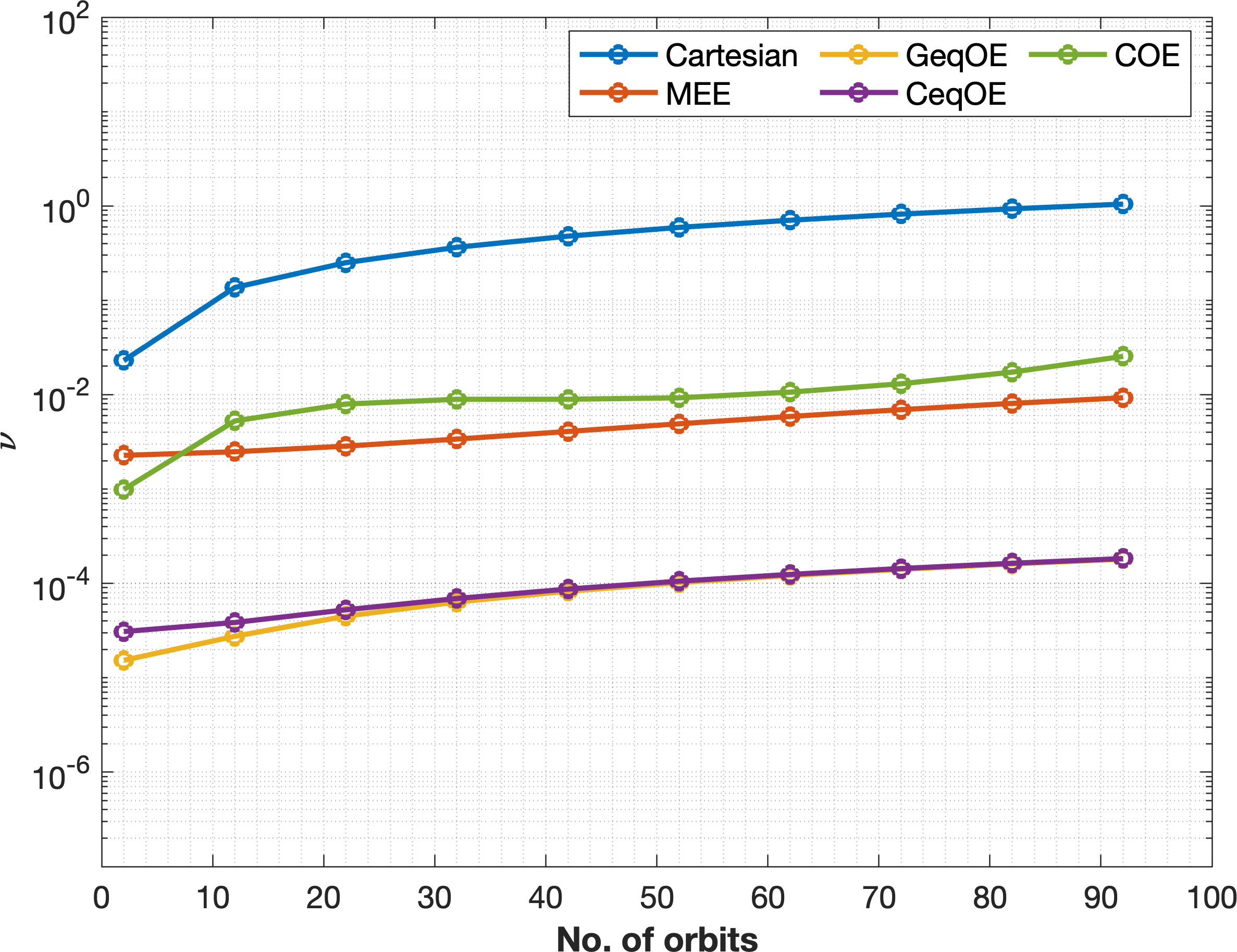}
    \caption{Comparison of the nonlinearity index (as defined in \cite{Junkins2003HowNI})  of GEqOE with other coordinates. }
    \label{GEqfig}
\end{figure}

 
{
To investigate the linearity of different coordinate systems, the nonlinearity index  $v$  proposed in \cite{Junkins2003HowNI} is used, where}

\begin{equation}
    v \triangleq \sup _{i=1, \cdots, N} \frac{\left\|A\left(\mathbf{x}_{\mathrm{i}}\right)-A(\overline{\mathbf{x}})\right\|}{\|A(\overline{\mathbf{x}})\|}.
\end{equation}
{
 $A$, $\mathbf{x}$ and $\overline{\mathbf{x}}$ are the state transition matrix, the perturbed state, and the unperturbed state after a propagation, respectively.  
A propagation was done for a spacecraft at an altitude of 350 km with an inclination of 99.22 deg in Cartesian, MEE, Generalized Equinoctial Orbital Elements (GEqOE) \cite{claudio}, Classical Equinoctial Orbital Elements (CEqOE) \cite{claudio2} and COE coordinates. Initial normalized position and velocity uncertainties of 1 km and 1 m/s were then introduced to analyze how the uncertainty is propagated forward in each coordinate system.  The propagated uncertainty of each of the coordinates was estimated using $v$. The results are shown in Figure \ref{GEqfig}, which illustrates that the nonlinearity index of GEqOE is lower than other coordinate systems, especially for a lower number of orbits propagated. Hence, GEqOEs are chosen to represent the convex nodes. 
}

Now the problem must be convexified to be solved using convex optimization. The convexification is done around the initial guess obtained in Section \ref{initialg}.  Firstly, the objective function Eq. \eqref{obj} can be convexified as follows

\begin{equation} \label{convexobj}
    K = \sum_{i=1}^{N_{seg}} \widetilde{\mathscr{a}}(i) {dt}(i) +  \widetilde{\Delta v'} \ \text{where }  dt(i) = t(i+1) - t(i)
\end{equation}
The problem dynamics can be convexified as 
\begin{equation}
\bm{x}_{GEq} ({i+1}) = \bm{A} (i)({\hat{\bm{x}}_{GEq}},\hat{\bm{\mathscr{a}}})\bm{x} (i) + \bm{B} (i)({\hat{\bm{x}}_{GEq}},\hat{\bm{\mathscr{a}}}) {\bm{\mathscr{a}} (i)} +\hat{\bm{x}}_{GEq}(i+1) - \bm{A} (i)(\hat{\bm{x}}_{GEq},\hat{\bm{\mathscr{a}}})\hat{\bm{x}}_{GEq}(i) - \bm{B} (i)(\hat{\bm{x}}_{GEq},\hat{\bm{\mathscr{a}}}) {\hat{\bm{\mathscr{a}} }(i)} 
\label{convexdyn}
\end{equation}

where $\bm{A} (i)({\hat{\bm{x}}}_{GEq},\hat{\bm{\mathscr{a}}})$ and $\bm{B} (i)({\hat{\bm{x}}}_{GEq},\hat{\bm{\mathscr{a}}})$ are evaluated at the initial guess trajectory, $[{\hat{\bm{x}}}_{GEq},\hat{\bm{\mathscr{a}}}]$. The subscript $_{GEq}$ denotes the fact that the states are given in GEqOE coordinates in the convex formulation. $i = 1,....,N_{seg}$ represents the convex nodes.  Note that $\bm{\mathscr{a}}_{t}^{convex}$ is written simply as $\bm{\mathscr{a}}$ from Eq \eqref{convexdyn} onwards for simplicity, and that the initial guess $\bm{\hat{x}}$ has been converted from Cartesian to GEqOE.  \par 
To ensure that the convexified dynamics do not cause the state to drift far from the reference- which would make the convexification invalid-, boundary conditions must be enforced on the convexified state variables. In this work, a trust region was enforced around the reference trajectory based on the nonlinearity of the dynamics involved. The Taylor-based nonlinearity index $\bm{\nu}_{\mathscr{T}}$ given in \cite{Losacco2023ALA}  is used to calculate a trust region.  $\bm{\nu}_{\mathscr{T}}$  is defined based on the deviation of the first order terms in the Taylor expansion of the state transition matrix $\bm{A}(\hat{\bm{x}}_{GEq}, \hat{\mathscr{a}})$. The boundary constraint on the state is enforced as 
 \begin{equation}
\begin{gathered} \label{boundary}
   \bm{x}_{GEq,p}(i)  \leq \hat{\bm{x}}_{GEq,p}(i)  +  {\zeta}/{\bm{\nu}_{\mathscr{T},p}(i)}  \\
     \bm{x}_{GEq,p}(i)  \geq \hat{\bm{x}}_{GEq,p}(i)  -  {\zeta}/{\bm{\nu}_{\mathscr{T},p}(i)}  
\end{gathered}
\end{equation}

where $\zeta$ is a user-set scalar,  and ${\bm{\nu}_{\mathscr{T},p}(i)}$ is a 6 element vector defined for each node.  $i = 1,..., N_{seg}+1$ is the current node and $p = 1,...,6$ represents the states of the node $i$.

The initial state condition is already convex, and given in GEqOE elements, it becomes 
\begin{equation}\label{cinit}
    \bm{x}_{GEq} (1)= \bm{x}_{0,GEq}
\end{equation}
The thrust acceleration constraint can be transformed into a convex second-order cone constraint as 

\begin{equation}\label{caccel}
 \sqrt{\mathscr{a}_r(i) ^2 + \mathscr{a}_\theta(i) ^2 + \mathscr{a}_n(i) ^2}  \leq \widetilde{\mathscr{a}}(i) 
\end{equation}
where $\widetilde{\mathscr{a}}(i) $ represents the thrust acceleration magnitude for node $i$. 
The thrust boundary becomes 
\begin{equation} \label{caccelb}
0 \leq \widetilde{\mathscr{a}}(i)  \leq \eta^{req}(i) T_{max}/m(i) 
\end{equation}

The $\Delta v'$ constraint must also be written as a convex cone. However, in order to obtain $\Delta v_a$ $\Delta v_h$ and $\Delta v_k$ defined on mean MEEs in a convex manner, a convex formulation of the conversion of the osculating GEqOE end state to mean MEEs is needed beforehand. Note that this involves combining the conversion function from GEqOE to MEE given in \cite{claudio} with the osculating to mean conversion given in \cite{8}.  Denoting this combined conversion function as $f: GEq \rightarrow \overline{MEE}$, and the Jacobian of the conversion function as $J_{f}$, a convex representation of $ \Delta \bar{\text{\oe}} = [\Delta \bar{a}, \Delta \bar{h}, \Delta \bar{k}]$  can be obtained as follows. 
\begin{equation}
\label{coeconv}
\begin{aligned}
    \Delta \bar{\text{\oe}} &=  [\bar{a}^c - \bar{a}^p, \bar{h}^c - \bar{h}^p, \bar{k}^c - \bar{k}^p] = [\bar{\text{\oe}}^c - \bar{\text{\oe}}^p] =  f (\bm{x}_{GEq}(N_{seg}+1)) - \bar{\text{\oe}}^p \\ 
       &= f (\hat{\bm{x}}_{GEq}(N_{seg}+1)) +  J_{f}(\hat{\bm{x}}_{GEq}(N_{seg}+1)) \bm{x}_{GEq}(N_{seg}+1) -  J_{f}(\hat{\bm{x}}_{GEq}(N_{seg}+1)) \hat{\bm{x}}_{GEq}(N_{seg}+1)    - \bar{\text{\oe}}^p 
\end{aligned}
\end{equation}
Following the convexification of the coordinate conversion, Eq \eqref{dvp2} is used to map $\Delta \text{\oe}$ to  $\Delta v_a$, $\Delta v_h$ and  $\Delta v_k$ variables. This equation is already convex in nature. Finally, the $\Delta v'$ constraint shall also be written as a convex cone

\begin{equation}
     \sqrt{\Delta v_a^2 + \Delta v_h^2 + \Delta v_k^2}  \leq \widetilde{\Delta v'}. \\
     \label{conedv}
\end{equation}
where $\widetilde{\Delta v'}$ represents the $\bm{\Delta v}'$ magnitude for node $i$. 

To summarize the convex problem following the convexification: 
     \begin{subequations}
            \begin{align}
                \mathop{\text{minimize}}_{\mathscr{a}(i)} ~~~~ & \text{Eq. \eqref{convexobj}} & \text{Minimization of control acceleration and $\Delta v'$.}\\
                \text{subject to}  ~~~~ &  \text{Eq. \eqref{convexdyn} for }  i = 1,...,N_{seg} &\text{Dynamics constraints in GEqOE}\\
                & \text{Eq. \eqref{boundary} for }  i = 1,...,N_{seg}+1 & \text{Trust region bound on states}\\   
                & \text{Eq. \eqref{cinit}}  & \text{Bound on the initial state} \\
                & \text{Eq. \eqref{caccel} for }  i = 1,...,N_{seg}   & \text{Second order cones for control acceleration}  \\
                & \text{Eq. \eqref{caccelb} for }  i = 1,...,N_{seg} &\text{Control acceleration bounds} \\
                & \text{Eq. \eqref{coeconv}}  &\text{Calculation of $\Delta \bar{\text{\oe}}$ from the end state} \\ 
                &\text{Eq. \eqref{dvp2}}  &\text{Calculation of the elements of $\bm{\Delta v'}$ from $\Delta \bar{\text{\oe}}$. } \\  
                & \text{Eq. \eqref{conedv}} &\text{ Second order cone for ${\Delta v'} = |\bm{\Delta v'}|$ } 
                \end{align}
            \label{eq:optiProb}
        \end{subequations}

Following the convexification, a convex solver such as MOSEK \cite{mosek} or Gurobi \cite{gurobi} can be utilized to solve it and obtain the optimized RTN thrust acceleration profile, denoted $\bm{\mathscr{a}}_t^{Convex}$. In this work, the state transition matrices involved were obtained using the Differential Algebra Computational Toolbox (DACE) \cite{dace} in C$++$, and MOSEK was used to perform convex optimization in Matlab. \par 
 The following outputs from the convex formulation are necessary for other steps of the convex-MPC guidance: the convex optimized acceleration for the segment $\bm{\mathscr{a}}_t^{Convex}$ and the PMDT reference orbital elements at $t^{convex}(N_{run} + N_{seg}+1)$, $[\bar{a}^{p}, \bar{h}^{p},\bar{k}^{p}]$.
\subsubsection{Forward Propagation of the Spacecraft} \label{fprops}

Once the optimal acceleration profile is obtained using convex optimization, the spacecraft is propagated under nonconvex and realistic dynamics and thrust uncertainties. In this section, thrust bias and magnitude errors are considered, as well as events of misthrust.\par 
Firstly, at the initialization stage, a misthrust profile vector  $\bm{M} = M_1, M_2, M_3, ..., M_m$ with one value per convex segment is generated by drawing from a uniform distribution between 0 and 1: 

\begin{equation}
 \bm{M}_i \sim \mathcal{U}(0,1) \ \text{where} \ i = 1:m
\end{equation}

 Note that $\bm{M}$ will therefore have $m = \text{ceil}(\text{len}(T)/N_{seg})$ values between 0 and 1. If $\bm{M} < p_{misthrust}$, where $p_{misthrust}$ is a user-set probability, the spacecraft is unable to thrust during the segment. A parameter $d_{misthrust}$ is also introduced, which denotes the user-set length of a misthrust event. Hence, if segment $i$ encounters a misthrust, it is extended till segment $i+ d_{misthrust}-1$. \par 
Algorithm \ref{alg:6} details the forward propagation process for a single segment.  Inputs for this algorithm are the initial osculating coordinates $\bm{x}_{0,CC}$, starting mass $m_0$, the convex-optimized RTN thrust acceleration profile $\mathscr{a}_t^{convex}$, PMDT reference at the segment end time $[\bar{a}^{p}, \bar{h}^{p},\bar{k}^{p}]$ and the misthrust fraction $M$.  {The fractional thrust magnitude error $\delta_T$ and the out-of-plane thrust angle error $\delta_\beta$ are drawn from normal distributions of user-set standard deviations $\bar{\delta}_T$ and $\bar{\delta}_\beta$, respectively.} Note that a recomputation of the eclipse profile $\eta^{req}$ is conducted in Algorithm \ref{alg:6}.  It does not reuse the thrust profile calculated in the initial guess generation in Algorithm \ref{alg:2}. This ensures that the thrust profile remains realistic and adheres to the eclipse-based criteria even in the presence of thrust errors and long-duration misthrust events that deviate the real trajectory from the reference. \par 
Also note that in Alogirhm \ref{alg:6}, the spacecraft is separately propagated under low and high-fidelity dynamical models. The low-fidelity model has the same dynamics as the convex optimization and considers only the accelerations due to drag, $J_2$ perturbations, and thrust. The low fidelity model - as well as the PMDT and the convex tracking- operates in the TOD frame, where the effects of precession and nutation are included.  The high-fidelity model used is the one developed in \cite{MORSELLI2014490}, which takes the full geopotential (zonal, tesseral and sectorial harmonics)  {up to order 8}, solar radiation pressure and the third body perturbations of the sun and the moon into account. and operates in the J2000 frame. Thus, a conversion from TOD to J2000 frame is implemented before the high-fidelity propagation, and a conversion from J2000 to TOD is conducted following the propagation.

\begin{algorithm}[hbt!]
\caption{Forward propagation for a convex tracking segment}\label{alg:6}
\textbf{Input} $\bm{x}_{s,CC} =  [\bm{r}_s , \bm{v}_s]^T = \boldsymbol{x}_{0,CC} $, $m_0$, $\bm{\mathscr{a}}_{t}^{convex}$, $M$, $[\bar{a}^{p}, \bar{h}^{p},\bar{k}^{p}]$
\begin{algorithmic}
\State Let $m(1) = m_0$. 
\If{$M < p_{misthrust}$} 
$ \bm{\mathscr{a}}_{t}^{convex}= \bm{0}$  \Comment{$ p_{misthrust}$ is the user-set misthrust probability}
\EndIf
\State Calculate the mean COE coordinates.  $\bar{\bm{x}}_{s,COE} = oscCart2meanKep(\boldsymbol{x}_{s,CC})$. \Comment{conversion given in \cite{8}}
\State  Calculate the mean argument of latitude: $\bar{L}_{fp}(1) =  \bar{\bm{x}}_{s,COE} (5) + \bar{\bm{x}}_{s,COE}(6) $
\State Define the time vector for the segment. $t = t^{convex}(N_{run} : N_{run} + N_{seg}) $ and let $m(1) = m_0$. 
\For{$i = 1: N_{seg}$} 
\State Compute the real eclipse/Duty cycle $\eta^{req}$ using $DC^{req}$ as given in Algorithm \ref{alg:3}.
\State Introduce normalized thrust magnitude and direction errors.
\State \hspace{\algorithmicindent} Compute the thrust magnitude with error:
\State \hspace{\algorithmicindent} \hspace{\algorithmicindent}  $  |\bm{\mathscr{a}} ^{convex}_{t} (i)|_{E}   = |\bm{\mathscr{a}} ^{convex}_{t} (i)|   + |\bm{\mathscr{a}} ^{convex}_{t} (i) | \delta_T \ \text{where}  \ \delta_T  \sim\mathcal{N}(0,\bar{\delta_T})$  \State\hspace{\algorithmicindent} \hspace{\algorithmicindent}  \Comment{$\bar{\delta_T}$:fractional standard deviation of thrust magnitude error}
\State \hspace{\algorithmicindent} Compute the out-of-plane angle:  $\beta = \sin^{-1}( {\bm{\mathscr{a}} ^{convex}_{t}(i)_n}/{|\bm{\mathscr{a}} ^{convex}_{t} (i)|})$
\State \hspace{\algorithmicindent} Compute the in-plane angle:  $\alpha = \tan^{-1}( {\bm{\mathscr{a}} ^{convex}_{t}(i)_\theta}/{\bm{\mathscr{a}} ^{convex}_{t}(i)_n})$
\State \hspace{\algorithmicindent} Calculate the out-of-plane angle with error: 
\State \hspace{\algorithmicindent}  \hspace{\algorithmicindent}  $\beta_E    = \beta  + \delta_\beta   \ \text{where}  \  \delta_\beta \sim \mathcal{N}(0,\bar{\delta_\beta})$ \Comment{$\bar{\delta_\beta}$: standard deviation of the out of plane thrust angle error}
\State \hspace{\algorithmicindent} Compute the thrust acceleration with errors: 
\State  \hspace{\algorithmicindent}  \hspace{\algorithmicindent}  $\bm{\mathscr{a}} ^{convex}_{t} (i)_E = |\bm{\mathscr{a}} ^{convex}_{t} (i)|_{E}  [\cos{\beta_E}\sin{\alpha}, \sin{\alpha}\cos{\beta_E}, \sin{\beta_E}]^T$
\State \hspace{\algorithmicindent} {Calculate $\Delta v = \eta^{req} \bm{\mathscr{a}} ^{convex}_{t} (i)_E(t(i+1) - t(i))$}
\State \hspace{\algorithmicindent}  Update mass. $ m(i) = {m(i-1)}/{\exp(\Delta v/I_{sp}/g_0)}$.
\If {Low fidelity propagation}
\State  {Propagate $\bm{x}_{s,CC}$ from $t(i \rightarrow i+1)$ under $\bm{a}_g , \bm{a}_{J_2} ,\bm{a}_d $ and $\eta^{req}  g(\bm{a}_{t}^{convex}(i))$ to get the state at $t(i+1)$.}  \Comment{As before $g:RTN \rightarrow TOD$.}
\EndIf
\If {High fidelity propagation}
\State Convert $\bm{x}_{s,CC}$ from TOD frame to J2000. 
\State  {Propagate $\bm{x}_{s,CC}$ from $t(i \rightarrow i+1)$ under $ \eta^{req} g(\bm{a}_{t}^{convex}(i))$ and the high fidelity dynamics model \\ \pushcode[1] in \cite{MORSELLI2014490}.}  
\State Convert the output state from J2000 to TOD.
\EndIf
\State  Define the new cartesian state at t(i+1) as $\bm{x}_{s,CC}$.
\State Calculate the mean COE coordinates.  $\bar{\bm{x}}_{s,COE} = oscCart2meanKep(\boldsymbol{x}_{s,CC})$. \Comment{conversion given in \cite{8}}
\State  Calculate the mean argument of latitude: $\bar{L}_{fp}(i+1) =  \bar{\bm{x}}_{s,COE} (5) + \bar{\bm{x}}_{s,COE}(6) $
\EndFor
\State Convert $\bar{\bm{x}}_{s,COE}$ to mean MEE. $\bar{\bm{x}}_{s,MEE}  = Cart2MEE(\bar{\bm{x}}_{s,COE})$
\State Calculate $\Delta v'_{fp}$, the $\Delta v'$ to go from $\bar{\bm{x}}_{s,MEE}$ to the target state for this segment $[\bar{a}^{p}, \bar{h}^{p},\bar{k}^{p}]$ using Eq. \eqref{dvp2}.
\end{algorithmic}
\textbf{Output } $\boldsymbol{x}_{f,CC} = \bm{x}_{s,CC}, \Delta v'_{fp}$ and $m_f = m(N_{seg})$, $\bar{L}_{fp}$.
\end{algorithm}

The following outputs from the forward propagation are necessary for other steps of the convex-MPC guidance: the end state of the propagation $\boldsymbol{x}_{f,CC}$, the $\Delta v'$ to go to the target state of this segment $\Delta v'_{fp}$, end mass $m_f$ and the mean true longitude profile  $\bar{L}_{fp}$.

\subsubsection{Phase Matching by Shifting the Segment Endpoints in Time} \label{Ltrack}

Following the forward propagation of an up leg, the last node of the convex segment is shifted in time such that the mean phasing angle (true longitude) can match up with the desired profile $\bar{L} ^{PMDT}_{adj}$, calculated in Section \ref{Ladjcal}.\par 
Firstly, the required time adjustment $dt$ that corresponds to shifting the true longitude by $d \bar{L}$ is calculated as follows, where $d \bar{L} = \bar{L}_{fp}(N_{seg}+1) - \bar{L}^{PMDT}_{adj}(N_{run} + N_{seg}+1)$ and $\dot{L}$ is the time derivative of the true longitude under natural dynamics.

\begin{equation}
    dt = k \frac{\cos^{-1}(\cos{d \bar{L}})}{\dot{L}}, \ \text{where} \ k = \text{sgn}(\sin(d \bar{L})/(1-\cos(d \bar{L})^2)^{1/2}
\end{equation}

$k$ is used to determine whether if the endpoint of the tracking segment needs to be shifted forwards or backwards in time. $\cos^{-1}(\cos{d \bar{L}})$ term is introduced to obtain the positive acute angle from $d \bar{L}$.
Now, the last node of the tracking segment is shifted from time $t(N_{seg}+1)$ to $t(N_{seg}+1) + dt$.  The spacecraft is propagated (without thrust acceleration, in osculating elements) from  $t(N_{seg})$ to $t(N_{seg}+1)+dt$ to obtain the new end state, $\bm{x}_{f, adj, CC}$. \par 
The following output from the phase matching is necessary for other steps of the convex-MPC guidance: adjusted end state in osculating, Cartesian coordinates $\bm{x}_{f,adj, CC}$.

\subsubsection{Computation of Recovery Time} \label{recT}

As the error builds up over time due to the presence of thrust uncertainties and nonlinear dynamics, the  {$\Delta v'_{fp}$} value calculated in Section \ref{fprops} increases.  However, most increases in $\Delta v'_{fp}$ can be brought back to a nominal value within a few tracking segments,  as the reference duty cycle $DC^{PMDT}$ being lower than the spacecraft's real duty cycle $DC^{req}$ provides a margin of thrust for error recovery.\par 

\begin{figure}[hbt!]
    \centering
 \includegraphics[width = 0.5\textwidth]{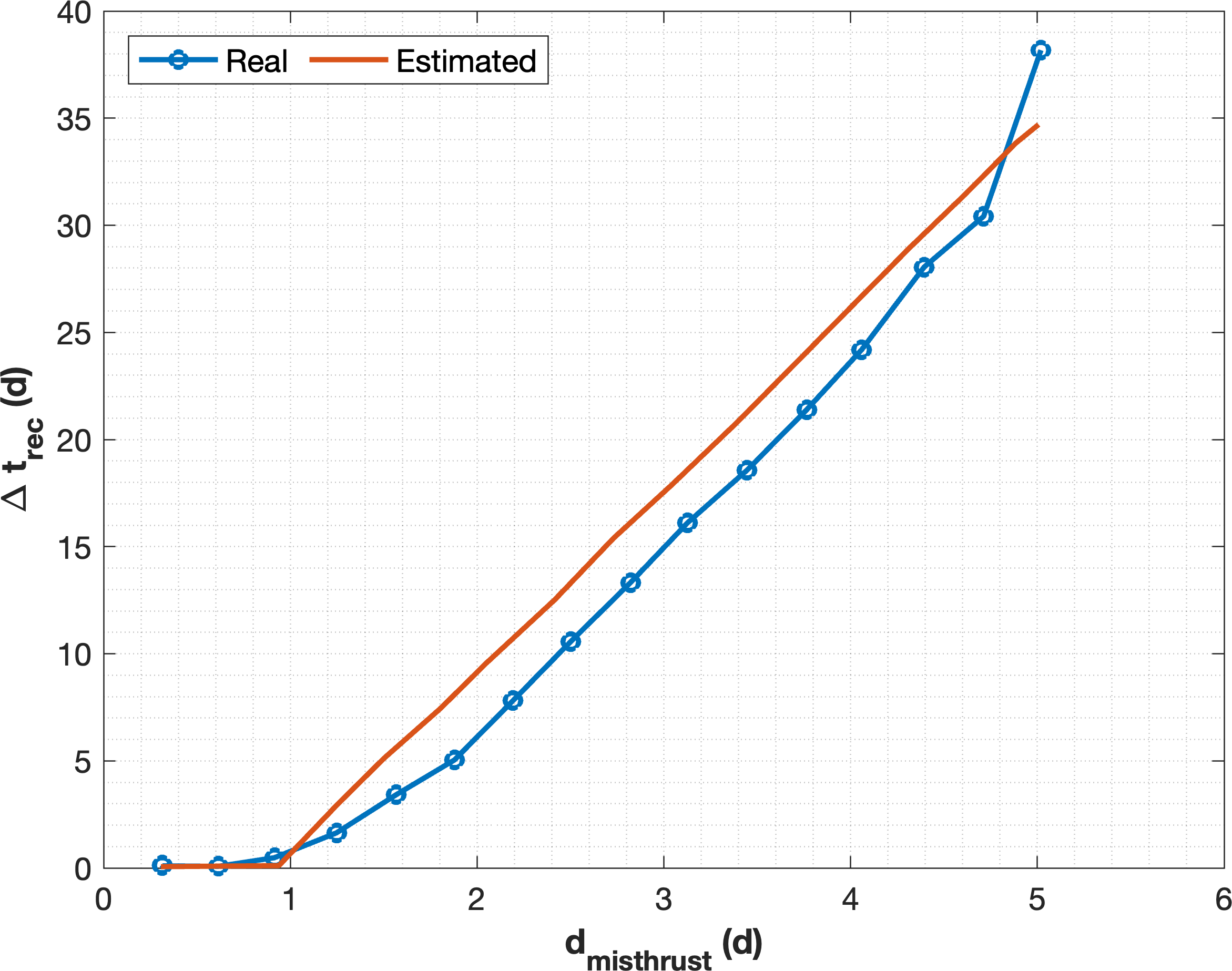}
    \caption{ Comparison of the time taken for misthrust recovery (real vs estimated using Eq. \eqref{dTrec})}
    \label{misthrust}
\end{figure}

The time taken for bringing the  $\Delta v'_{fp}$ back to a user-set $\Delta v'$ threshold- denoted as $\Delta v'_T$- can be calculated as a function of the thrust margin as  
\begin{equation}\label{dTrec}
    \Delta t_{rec} = \frac{\Delta v'_{fp} - \Delta v'_{T}}{ (DC^{req} - DC^{PMDT}) T_{max}/m_f}
\end{equation} 

where $m_f$ is the mass remaining at the end of the segment, as defined in Section \ref{fprops}. To determine if Eq. \eqref{dTrec} is a good estimate of the recovery time, the real recovery times for different misthrust durations were compared against the times estimated by Eq. \eqref{dTrec}. For this comparison, $\Delta v'_T$ was set to be $\SI{2}{\meter\per\second}$, and $DC^{req}$ and $DC^{PMDT}$ were set to 0.50 and 0.40, respectively. The results are given in Figure \ref{misthrust}, which shows that Eq. \eqref{dTrec} is a close upper estimate for the real recovery time for most misthrust durations. Hence, it is likely a reliable parameter to gauge whether a recomputation is needed to reach the target.\par 
The following outputs from the recovery time calculation are necessary for other steps of the convex-MPC guidance:  the recovery time $\Delta t_{rec}$ and the remaining time $t^{convex}_{rem} = t_f - t^{convex}(N_{run} + N_{seg} +1)$.

\subsubsection{Accuracy Check}\label{Acheck}
If $\Delta t_{rec} > t^{convex}_{rem}$- which is likely if a long-duration misthrust event takes place towards the end of the mission- the algorithm would require a new reference to be generated, as discussed in Section \ref{RefRegen}. Otherwise, it would proceed to the following check.

\subsubsection{Last Node check}\label{Lcheck}
The algorithm has reached the last node if $N_{run} \geq len(t^{convex})$. In this case, it would stop. If the last node has not been reached, the algorithm will proceed to the next convex tracking segment by setting $\bm{x}_{0,CC} = \bm{x}_{f,CC}$ (if phase matching is not conducted) or $\bm{x}_{0,CC} = \bm{x}_{f,adj,CC}$ (if phase matching is conducted), $m_0 = m_f$ and $N_{run} = N_{run} +N_{seg}$ and going back to Section \ref{initialg}. 

\subsection{Reference Recomputation} \label{RefRegen}

\begin{figure}[hbt!]
    \centering
    \includegraphics[width = 0.8\textwidth]{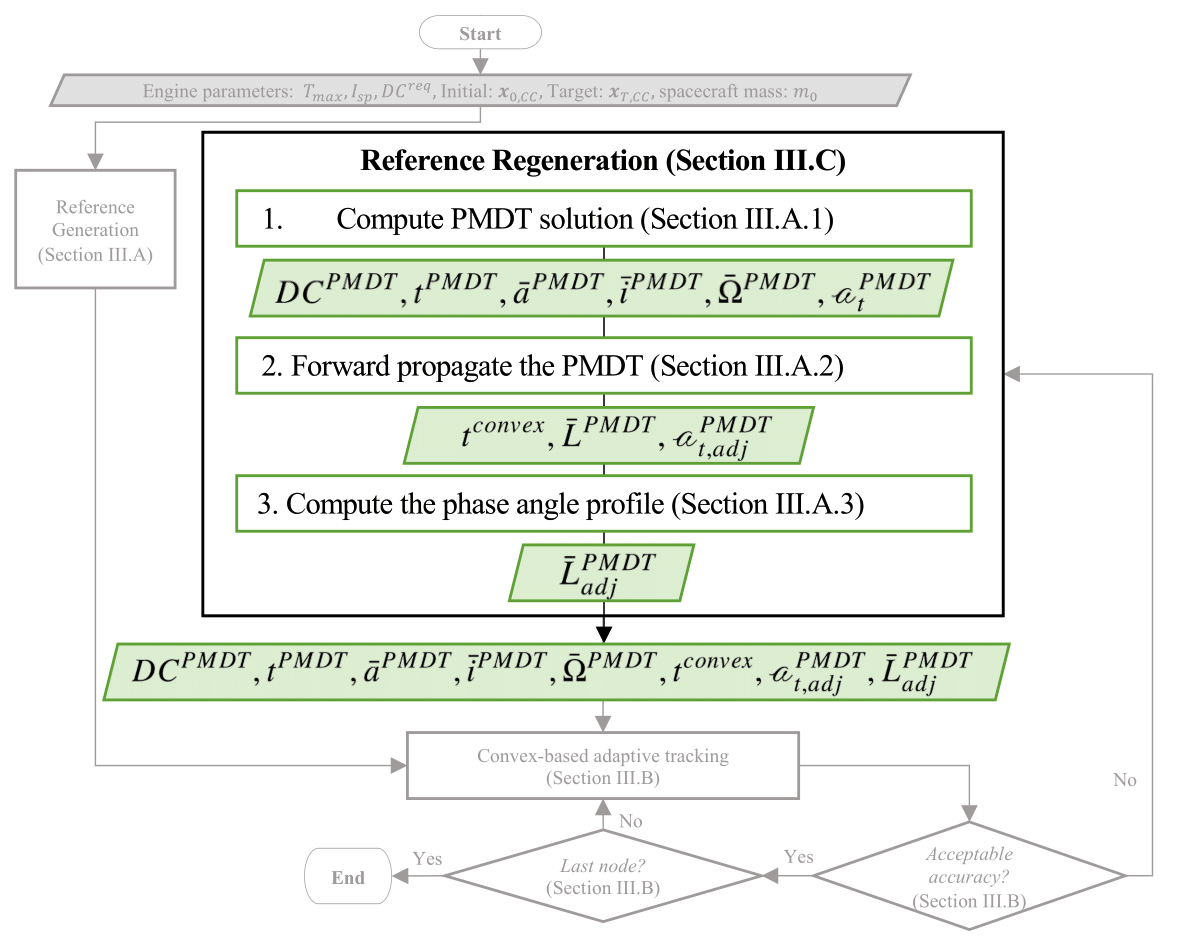}
    \caption{Reference Regeneration process}
    \label{RRgen}
\end{figure}

If a recomputation is requested (i.e., if $\Delta t_{rec} > t^{convex}_{rem}$), the PMDT reference and its forward propagation must be recomputed to obtain a new trajectory that allows the spacecraft to reach the target at an incremented flight time and/or $\Delta v$. All the steps conducted in Section \ref{Refgen} must be repeated when regenerating a new reference, as shown in Figure \ref{RRgen}. Note that if a recomputation is needed for an up leg that requires phase matching, a new $DC^{PMDT}$ must also be chosen in the same manner as before to ensure that the difference between the targetted mean true longitude $\bar{L}_T (t_f) $ and the predicted mean true longitude from the PMDT forward propagation $\bar{L} ^{PMDT}(t_f)$ is sufficiently close to zero (i.e., $\delta \bar{L}(t_f) \approx 0 $). \par

\section{Numerical Simulations and Analysis} \label{results}

In this section, results are provided for an up and down leg of the multi-debris removal tour. The up leg takes the servicer spacecraft from a 350 km altitude orbit to the first aim point to reach an H-2A (F15)\footnote{Tracking data available at \url{https://www.n2yo.com/satellite/?s=33500}} rocket body.  The first aim point is defined as the end point of phasing and the starting point of the relative navigation operations of rendezvous, typically located a few kilometers below and a few tens of kilometers behind the target \cite{fehse}.  The down leg brings the debris and the servicer back down to a 350 km altitude orbit.  {Note that the up-leg consists of phase matching, as rendezvous with the target is its expected outcome. Down leg does not require phase matching, as its target is only to reach a 350-km altitude quasi-circular orbit.} \par 
 
The servicer spacecraft has a starting wet mass of $m_0 = 800$ kg, with a drag coefficient of 2.2 and a frontal area of $\SI{0.01}{\meter\squared}$. The low thrust engine of the Servicer has 60 mN maximum thrust, 1300 s $I_{sp}$, and a duty cycle $DC^{req}$ of 50\% during both legs. \par 
 
Once the spacecraft reaches the first aim point of the rocket body, thirty days are allocated to approach and dock with the target. Simulations were conducted under the thrust error conditions shown in Table \ref{tError}. Monte Carlo simulations  {with 100 samples} were conducted to analyze their effect. \par 

\begin{table}[hbt!]
\centering 
\caption{Thrust error settings utilized in the case study}
\label{tError}
\begin{tabular}{lcccc} \hline 
Simulation         & $d_{misthrust}$&  $p_{misthrust}$ (\%) &  $\bar{\delta_T}$ (\%) & $\bar{\delta_\beta}$(deg) \\ \hline 
 Perfect thrust   & - & 0  & 0   & 0  \\
 Low error & 1 segment & 1                 & 5                                                    & 5                                                                       \\
High error, short misthrust (SM) & 2 segments  & 1              & 7                                                  & 7   \\ 
High error, long misthrust (LM) & 4 segments &  1            & 7                                                 & 7 \\  \hline                                                         
\end{tabular}
\end{table}

In the convex tracking, the number of nodes per orbit N was set to 36. The number of orbits per segment (n) was set to 5. Given the nonlinearity index analysis from Figure \ref{GEqfig},  this was determined to be a good compromise between the efficiency and accuracy of the convex tracking using GEqOE. The trust region parameter $\zeta$ was set to be 0.01. The targeted accuracy of the guidance is set to be  $\Delta v'_T = \SI{1.97}{\meter\per\second}$, corresponding to a 0.05\% error in $a$, $i$ and $\Omega$.  This value determines if a recomputation is necessary via Eq \eqref{dTrec}. 

\subsection{Fuel Optimal Up Leg from the Initial Orbit to the Target Debris}

In this section, the up leg of the mission is optimized for fuel consumption using the PMDT, and the obtained trajectory is used for the convex-MPC guidance. Firstly, the location of the first aim point is determined. 

\subsubsection{First Aim Point Determination}
{
The process to calculate the first aim point is given in Algorithm \ref{Firstaim}. The inputs to  Algorithm \ref{Firstaim} are the TOD coordinates of the debris $\bm{x}_{deb,CC}$ and the placement of the first aim point with respect to the debris in RTN coordinates $\bm{\delta r}_{CC}$. In this study, the first aim point was chosen to be 3 km below and 100 km behind the target debris, in accordance with the criteria specified by Astroscale in \cite{gnc}. Hence $\bm{\delta r}_{CC}= [\SI{-3}{\kilo\meter}, \SI{-100}{\kilo\meter}, 0]^T$. The output of Algorithm \ref{Firstaim} are the cartesian coordinates of the first aim point $\bm{x}_{fa,CC}$ .}

\begin{algorithm}[hbt!]
 \textbf{Input} $\bm{x}_{deb,CC} = [\bm{r}_{deb,CC}, \bm{v}_{deb,CC}]$, $\bm{\delta r}_{CC}$
\caption{First Aim point calculation }\label{Firstaim}
\begin{algorithmic}

\State Compute the first aim point position vector  $\bm{r}_{FA,CC} = \bm{r}_{deb,CC} + g(\bm{\delta r}_{CC})$ where $g: RTN \rightarrow TOD$.
\State Compute the Classical Orbital Elements of the debris $\bm{x}_{deb,COE} =  cart2coe(\bm{x}_{deb,CC})$. \Comment{Where $\bm{x}_{deb,COE} = [a_{deb}, e_{deb}, i_{deb},\Omega_{deb}, \omega_{deb}, \theta_{deb}]. $}
\State Note that the First aim point is in a co-elliptic orbit of the same orientation as the debris. Hence:  $e_{fa} = e_{deb}, i_{fa} = i_{deb}, \Omega_{fa} = \Omega_{deb}$ and $ \omega_{fa} = \omega_{deb}$
\State Convert $\bm{r}_{FA,CC}$ to periforcal frame to generate $\bm{r}_{FA,CC,pf}$, using the conversions given in \cite{perifocal}.
\State Calculate the true anomaly of the first aim point $\theta_{fa} = \tan^{-1}(\frac{\bm{r}_{FA,CC,pf}(2)}{\bm{r}_{FA,CC,pf}(1)})$.
\State Calculate the semi-major axis of the first aim point $a_{fa} = |\bm{r}_{FA,CC,pf}| ({1 + e_{fa} \cos{\theta_{fa}}})/({1 - e_{fa}^2})$.
\State Obtain the first aim point velocity in perifocal frame $\bm{v}_{fa,CC,pf} = [-\sqrt{{\mu}/p_{fa}} \sin{\theta_{fa}}, \sqrt{{\mu}/(p_{fa})} (e_{fa} + \cos{\theta_{fa}}),0]^T$.  \Comment{Note $ p_{fa} = a_{fa}(1- e_{fa}^2)$.}
\State Convert the velocity to TOD, $\bm{v}_{fa,CC} = g^{-1} (\bm{v}_{fa,CC,pf})$. 
\end{algorithmic}
\textbf{Output} $\bm{x}_{fa,CC} = [\bm{r}_{fa,CC} ;\bm{v}_{fa,CC}]$. \Comment{Note $\bm{x}_{T,CC} = \bm{x}_{fa,CC}$. }
\end{algorithm}

\subsubsection{Input Parameters}

The up leg starts at 00:00 UTC on 25 March 2022. The mean orbital parameters of the  H-2A (F15) rocket body{, the first aim point} and the starting orbit at this time are given in Table \ref{T1}. Semi-major axes, inclination, RAAN, and phasing angle ($\bar{L}=  \bar{\Omega} + \bar{\omega}+ \bar{\theta}$) are provided as targets, and the eccentricity is maintained to be quasi-circular by reversing the thrust direction across the line of nodes, as mentioned in Algorithm \ref{alg:3}.

\begin{table}[hbt!] 
\caption{Orbital elements of  {initial state, debris state, and first aim point state at the start time }}
\label{T1} 
\centering\begin{tabular}{lcccccc} \hline Object & $\bar{a}$ (km) & $\bar{e}$ & $\bar{i}$ (deg) & $\bar{\Omega}$ (deg) & $\bar{\omega}$ (deg) & $\bar{\theta}$ (deg) \\ \hline Initial & 6718.436 & 0.00348 &  98.306 &  15.200 &  0.0 & 0.0 \\ 
H-2A (F15) debris & 6980.031 &  0.00473 & 98.109 &  20.260  & 185.310 & 293.711 \\
First aim point & 6977.577 &  0.00473 & 98.109 &  20.260  & 184.823 & 293.374
\\\hline 
\end{tabular}
\end{table}

A duty cycle $DC^{PMDT}$ of 0.4019 was selected for the reference generation, as it provided a sufficient thrust margin while also generating a trajectory that closely approaches the true longitude target.  With this $DC^{PMDT}$ , the PMDT generates an optimal trajectory with $\Delta v = \SI{150.178}{\meter\per\second}$ and $TOF = \SI{59.182}{}$ days for the up leg.  A drift period of 1.853 days is needed to reach the target RAAN. At the end of the transfer, the forward propagated true longitude is only 8.31 deg away from the target.  Note that if $DC^{PMDT}= 0.5$, the resultant $\Delta v$ and $TOF$ are $\SI{150.078}{\meter\per\second}$ and $TOF = \SI{59.436}{}$ days, respectively. While the higher $DC^{PMDT}$ will shorten the thrust arcs, the drift time is increased to 13.356 days to match the target RAAN. Furthermore, the forward propagated true longitude is  112.57 deg away from the target, which is above the required 45 deg margin. Hence, if 0.5 were chosen instead as the reference duty cycle, phase matching would be harder and the reference would be made invalid relatively soon after departure due to the large shifts in the convex-tracking endpoints.  \par

\subsubsection{Up Leg Results}

\begin{figure}[hbt!]
  \centering
  \begin{subfigure}[b]{0.49\textwidth}
  \centering
    \includegraphics[width = \textwidth]{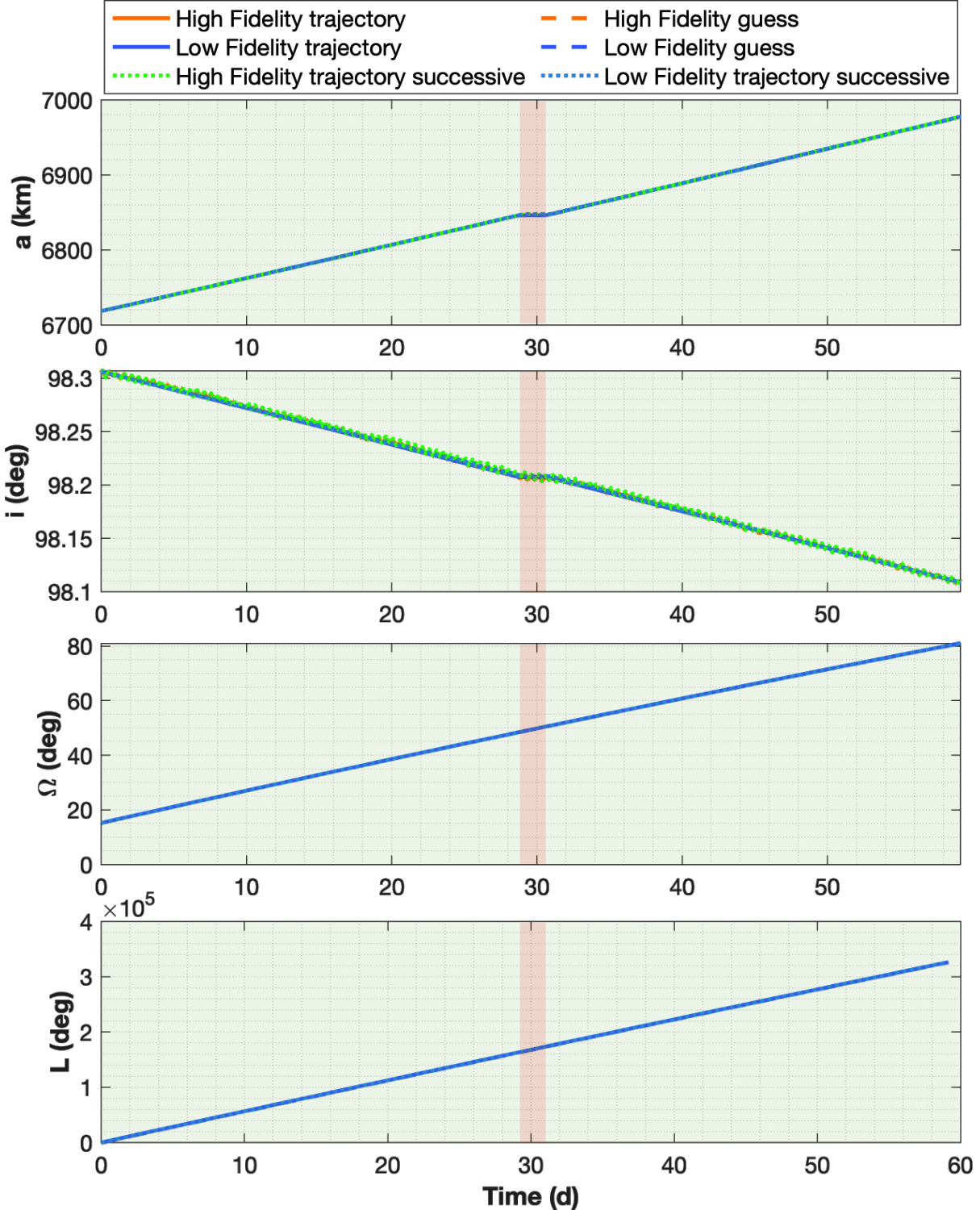}
    \caption{The fuel optimal solution with the high and low fidelity spacecraft propagations}
    \label{foptb1}
  \end{subfigure}
  \hfill
  \begin{subfigure}[b]{0.49\textwidth}
    \centering
    \includegraphics[width = \textwidth]{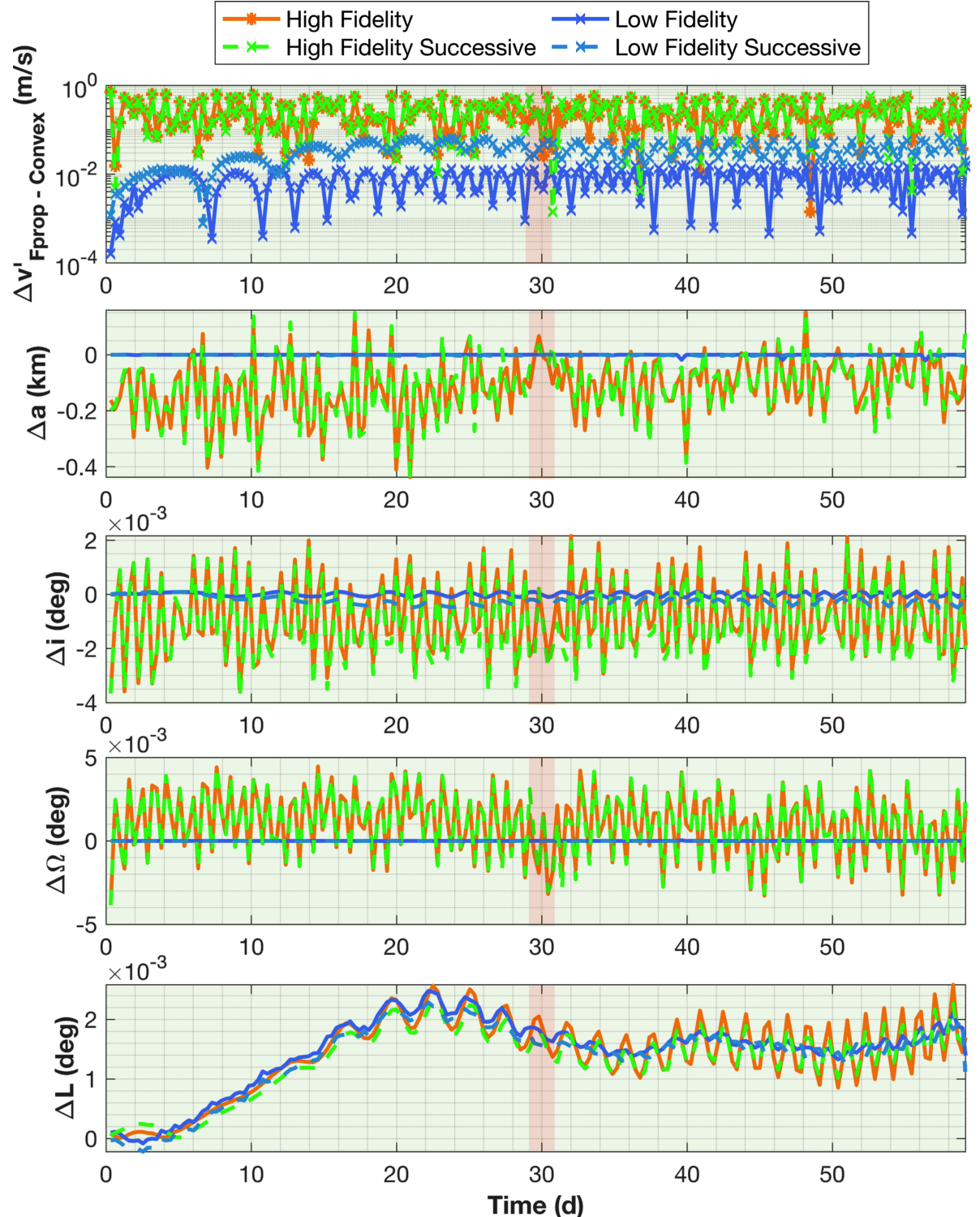}
    \caption{The evolution of $\Delta v'$, semi-major axis, inclination, and RAAN errors with respect to the reference trajectories over time}
    \label{foptb2}
  \end{subfigure}
  \caption{Up leg fuel-optimal results (perfect thrust) (Thrust arcs are highlighted in green, and the intermediate drift arc is given in red.)}
  \label{fig:main1}
\end{figure}

Figures \ref{foptb1} and \ref{foptb2} show the trajectory and error profiles encountered under perfect thrust, which illustrate that good semi-major axis, inclination, RAAN, and phasing angle accuracy can be maintained throughout the trajectory with both dynamics. The 1.853-day drift period is indicated by the flat region between the two thrust arcs in Figure \ref{foptb1}. The low-fidelity trajectory experiences less deviation from the reference than the high-fidelity one, as the latter encounters more deviations due to high-order perturbations, third-body effects, and solar radiation pressure. \par 

{Figures \ref{foptb1} and \ref{foptb2} also compare the outcomes of single iteration vs successive iteration convex tracking. In the case of successive iterations, the steps in section \ref{convex} followed by section \ref{fprops} are run in a loop, till the forward propagated trajectory and the convex optimization show strong agreement- the loop was set to terminate when $\Delta v'_{fprop - convex} \leq \SI{0.001}{\meter\per\second}$ or exceeds five iterations.  It can be seen that the successive and single iteration cases are similar in terms of accuracy, as shown in the top Figure in  Figure \ref{foptb2}.  The results obtained are compared in Table \ref{tR}, which shows that successive convexification does result in a minute increase in accuracy, at the cost of a small $\Delta v$ increase, for both high and low fidelity cases. The time taken to execute one low-fidelity convex tracking segment is 2.94 seconds, while for a high-fidelity counterpart, it extends to 6.82 seconds for a single iteration. This duration notably increases to 8.38 seconds and 15.96 seconds, respectively, for a successive iteration segment. Note that  $\sim 85$ \% of the computational time is occupied by the initial guess generation and forward propagation processes.  Hence, the increment in accuracy is negligible compared to the approximate tripling of computational time observed when successive convexification is used. } \par

\begin{figure}[hbt!]
    \centering
    \includegraphics[width = \textwidth]{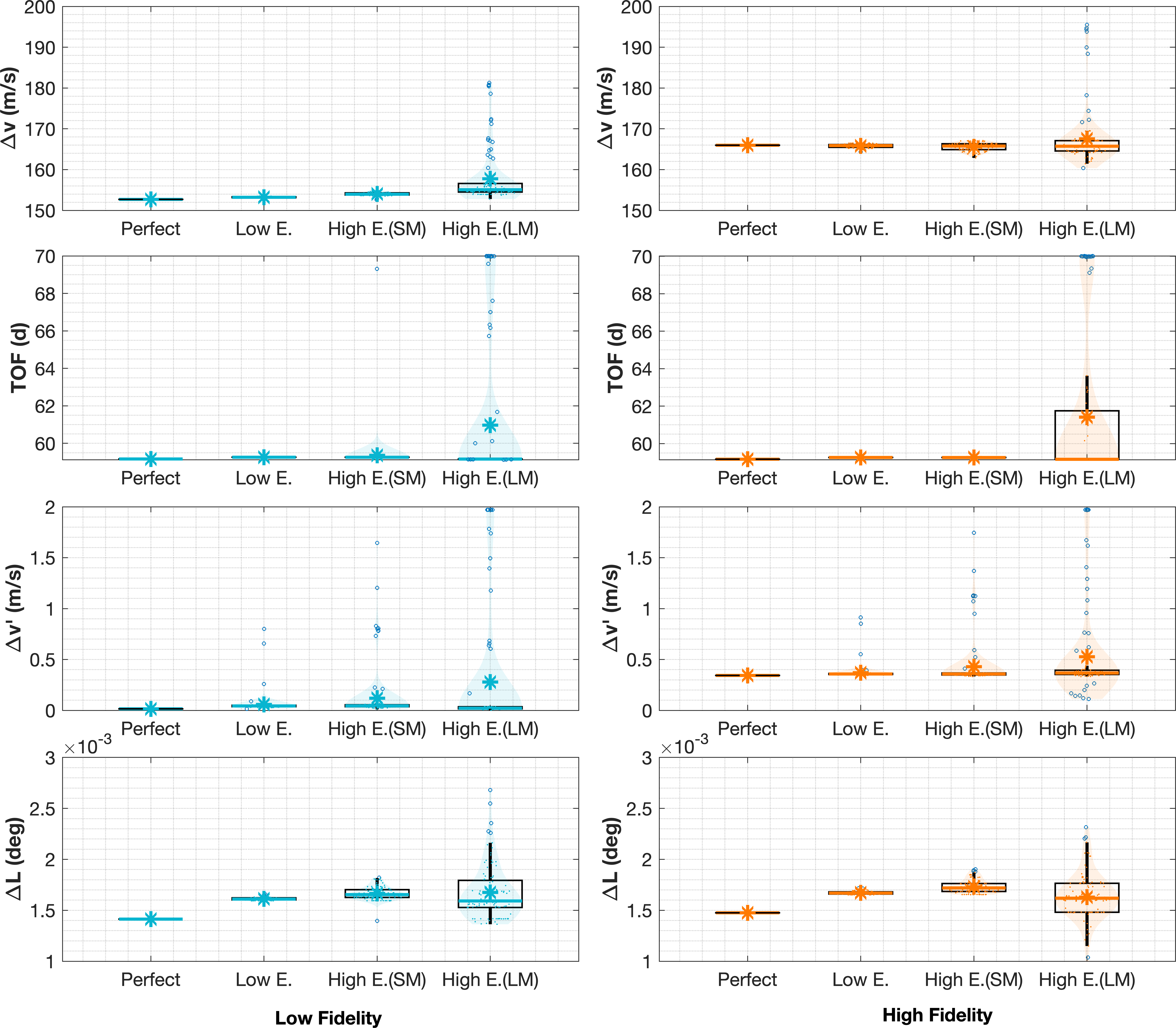}
    \caption{Up leg fuel-optimal results with thrust errors  (SM: short misthrust, LM: long misthrust) } 
    \label{foptb13}
\end{figure}

100 sample Monte Carlo simulations were conducted for cases with thrust errors, and the results are shown in Figure \ref{foptb13}.  {Note that the same realization of the 100 Monte-Carlo errors is sampled for each case with thrust errors.}

The numerical results of the perfect thrust case and the $25\%$ and $75\%$ quantiles of the Monte-Carlo simulations are detailed in Table \ref{tR}. Note that the $\Delta \bar{a}, \Delta \bar{i}, \Delta \bar{\Omega}$ and $\Delta \bar{L}$ are the errors with respect to the first aim point at the end of the transfer. \par 

As expected, due to the discrepancy of dynamics between the convex tracking and forward propagation of the spacecraft, the $\Delta v'_{fp}$s of the high-fidelity simulations are greater than that of the low-fidelity simulations.
From Table \ref{tR}, it can be seen that the fuel consumption of the high-fidelity simulations is greater than that of the low-fidelity simulations for all cases. This is due to two reasons: (1) More fuel is needed to adjust for deviations due to high-order perturbations, third-body effects, and solar radiation pressure in high-fidelity simulations. (2)   The PMDT and the low fidelity simulations use the Harris-Priester (HP) atmospheric density model \cite{HATTEN2017571} while the high fidelity propagation uses the more accurate NRLMSISE-00 \cite{MORSELLI2014490}. NRLMSISE-00 estimates higher densities than the HP model at low altitudes \cite{compare}. Hence, in high-fidelity dynamics, the higher atmospheric drag results in a larger fuel requirement to travel against it. \par  

It can be seen from Figure \ref{foptb13} that for each propagation model, increased thrust errors translate to increased upper limits of fuel consumption. This can be expected as propellant is required for error correction.  In the  {high error, long misthrust} cases, if the misthrusts occur towards the end of the transfer, recomputations are needed to reach the target. This is shown by the existence of the heavy tails and outliers in the $TOF$ and $\Delta v$ distributions shown in the high error, long misthrust cases in Figure \ref{foptb13}.   \par

\begin{table}[hbt!]
\centering
\caption{Up leg fuel-optimal results}
\label{tR}
\setlength{\tabcolsep}{1pt}
\adjustbox{max width=1.1\textwidth}{%
\begin{tabular}{cccccccc} 
\hline 
Spacecraft  & $\Delta \bar{a}$ & $\Delta \bar{i} $ & $\Delta \bar{\Omega}$ & $\Delta \bar{L}$ & TOF & $\Delta v$& $\Delta v^{'}_{fp}$ \\ 
Propagation & (km) & (deg) &  (deg) & (deg) &  (days) & (m/s) & (m/s) \\ \hline 
\multicolumn{8}{c}{PMDT}                                                                                                      \\
-& - & - & - &- & 59.182  & 150.18 & -\\ \hline 

\multicolumn{8}{c}{Perfect thrust  (Singe iteration) }      \\
Low Fidelity & -0.00152 & 2.47e-06 & -0.000496 & 0.00141 & 59.165 & 152.70 & 0.0152 \\
High Fidelity & -0.0358 &  0.00180 & -0.00215 & 0.00148& 59.164 & 165.97 & 0.343 \\\hline 
\multicolumn{8}{c}{Perfect thrust  (Successive iterations) }      \\
Low Fidelity & -0.00168 & -2.22e-05 & -0.000356 & 0.00133 & 59.140 & 152.78 & 0.0144 \\
High Fidelity & 0.0344 & 0.00235 & -0.00221 & 0.00137& 59.138 & 166.12 & 0.339 \\\hline
\multicolumn{8}{c}{Low error [$Q_1 \ Q_3$]}                                                               \\
Low Fidelity & [-0.0362 -0.00399] & [1.05e-05 0.000129] & [-0.000494 -0.000235] & [0.00143 0.00161]& [59.165 59.259] & [153.20 153.52] & [0.0162 0.0443] \\High Fidelity & [-0.0435 0.0529] & [-0.000792 0.00183] & [-0.00264 -0.00216] & [0.000154 0.00139]& [59.163 59.257] & [165.50 166.18] & [0.346 0.358] \\\hline

\multicolumn{8}{c}{High error (short misthrust (SM)) [$Q_1 \ Q_3$]}         \\
 Low Fidelity & [-0.0454 -0.0102] & [2.09e-05 0.00017] & [-0.000493 -0.000221] & [0.00148 0.00165]& [59.163 59.258] & [153.91 154.46] & [0.0194 0.0491] \\High Fidelity & [-0.0503 0.0453] & [-0.000719 0.0019] & [-0.00269 -0.00218] & [0.000153 0.0014]& [59.162 59.256] & [164.99 166.37] & [0.348 0.366] \\\hline

\multicolumn{8}{c}{High error (long misthrust (LM)) [$Q_1 \ Q_3$]}        \\
 Low Fidelity & [-0.0364 -0.00375] & [-3.06e-05 4.84e-05] & [-0.000552 -0.000435] & [0.00153 0.0018]& [59.157 59.165] & [154.50 156.65] & [0.0159 0.0359] \\High Fidelity & [-0.0548 -0.0317] & [0.000778 0.00204] & [-0.00225 -0.00212] & [0.00137 0.0014]& [59.158 61.669] & [164.57 167.09] & [0.353 0.394] \\\hline 
\end{tabular}}
\end{table}

\begin{figure}[hbt!]
  \centering
  \begin{subfigure}[b]{0.49\textwidth}
  \centering
    \includegraphics[width = \textwidth]{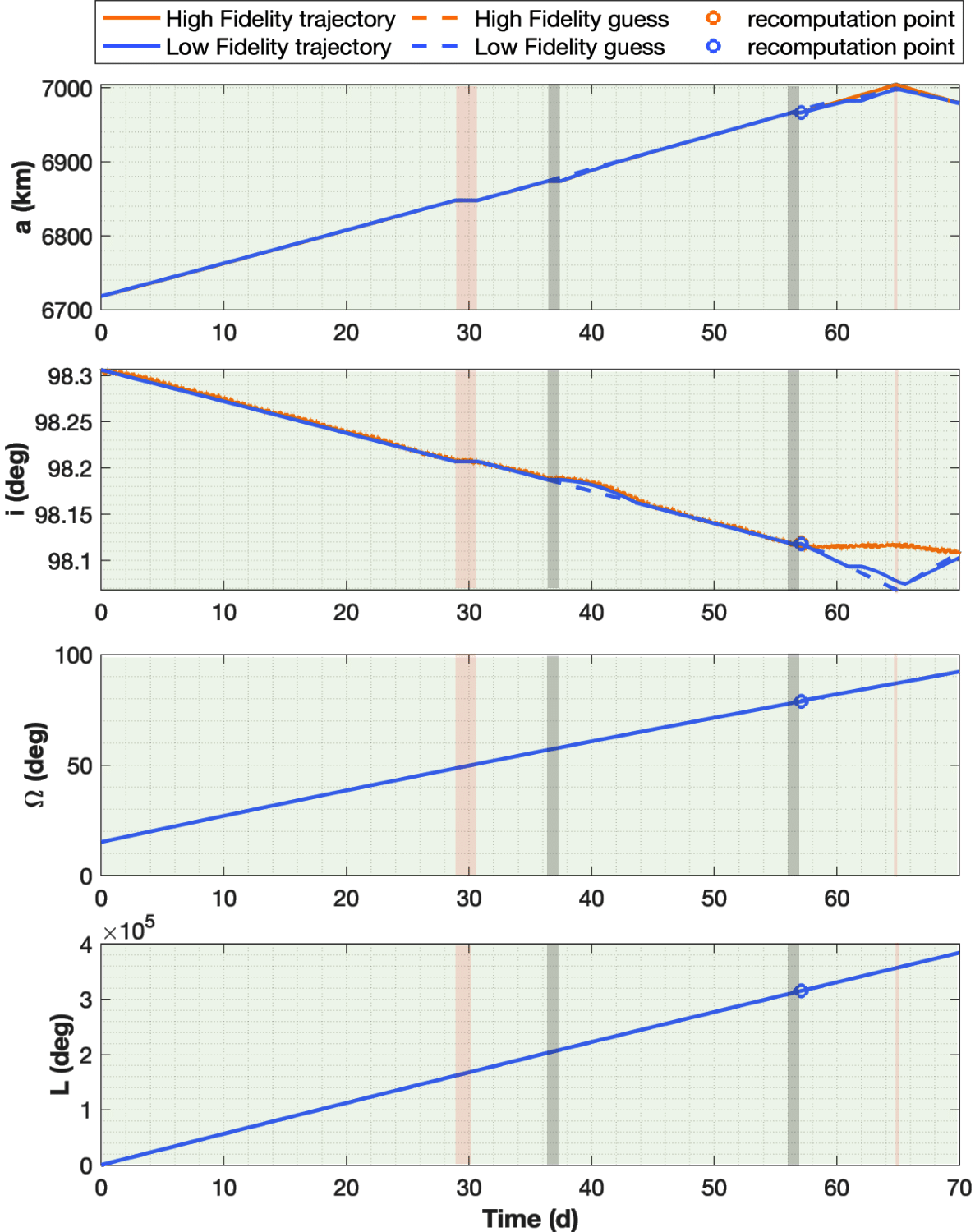}
    \caption{The fuel optimal solution with the high and low fidelity spacecraft propagations}
    \label{fopth1}
  \end{subfigure}
  \hfill
  \begin{subfigure}[b]{0.49\textwidth}
    \centering
    \includegraphics[width = \textwidth]{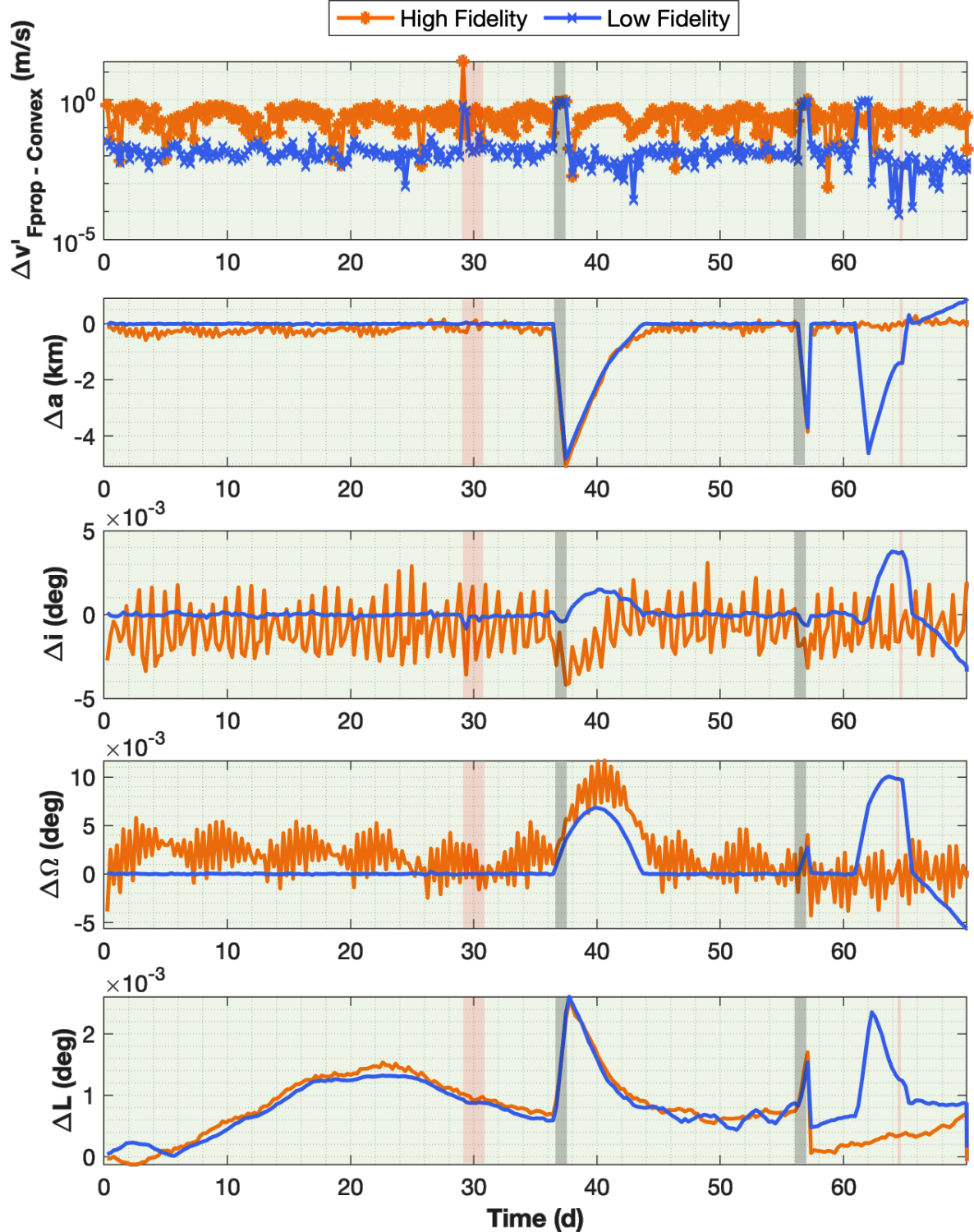}
    \caption{The evolution of $\Delta v'$, semi-major axis, inclination, and RAAN errors with respect to the reference trajectories over time}
    \label{fopth2}
  \end{subfigure}
  \caption{ {Up leg fuel-optimal results (high error LM). Thrust arcs are highlighted in green, and the intermediate drift arc is given in red. The long misthrusts are highlighted in grey.}}
  \label{fig:mainh}
\end{figure}

Figure \ref{fig:mainh} illustrates the trajectory and error profile of a high error LM simulation necessitating a recomputation to reach the target when the $\Delta v'_T$ is set to $\SI{1.97}{\meter\per\second}$. The initial trajectory experiences two LM events (highlighted in grey in Figure \ref{fig:mainh}. Both high and low-fidelity propagations are able to recover from the first misthrust event within a few convex tracking iterations. However, the second misthrust event happens when the spacecraft has travelled for 57.17 days of the 59.31-day initial transfer; hence, the time remaining is insufficient for natural recovery, and a recomputation is needed for both propagations. For the low-fidelity propagation, the new trajectory extended the $\Delta v$ and $TOF$ by $\SI{ 31.94}{\meter\per\second}$ and $\SI{11.023}{}$ days, respectively. For the high-fidelity propagation, the increments are  $\SI{ 46.14}{\meter\per\second}$ and $\SI{11.026}{}$ days. As expected, the high fidelity recomputation carries a greater added $\Delta v$ cost.   At the end of the transfer, the low-fidelity model yields a final $\Delta v'_{fp}$ of $\SI{0.916}{\meter\per\second}$, and the high-fidelity approach results in a $\Delta v'_{fp}$ of $\SI{1.841}{\meter\per\second}$.  This further underscores that while recomputations carry a considerable cost, they effectively ensure that the target is reached within the desired accuracy margins by the mission's end.

\subsection{Fuel Optimal Down Leg from the Target Debris to the Initial Orbit}

\subsubsection{Input Parameters}
At the start of the down leg, the total mass is the combination of the dry mass of the servicer, the remainder of the fuel mass \footnote{Note that the remaining fuel mass of the perfect thrust, low fidelity up leg was considered here for simplicity.}, and the debris mass ($m_0 = 790.38 +2991 = \SI{3781.4}{\kilo\gram}$). Thirty days are allocated for the rendezvous procedure; hence, the down leg start date is 00:00 UTC 20-Jun-2022. Table \ref{data2} shows the starting and target orbital parameters for the down leg. Note that during the down leg, only the semi-major axis is tracked and the eccentricity is maintained to be quasi-circular by alternating the direction of the out-of-plane thrust as shown in Algorithm \ref{alg:3}. \par 

\begin{table} [hbt!]\centering
\caption{Orbital elements of initial debris and the target orbit at launch time of the down leg}  \label{data2}   \adjustbox{max width=\textwidth}{%
\begin{tabular}{ccccccc} \hline  Object & $a$ (km) & e & $i$ (deg) & $\Omega$ (deg) & $\omega$ (deg) & $\theta$ (deg) \\ \hline Initial (H-2A F-15) & 6987.0507 & 0.0042309 & 98.2219 & 108.8944 & 275.8823 & 64.4907 \\Target & 6728.1363 & - & -& - & - & - \\\hline \end{tabular} }
\end{table}

A duty cycle $DC^{PMDT}$ of 0.45 was selected for the reference generation for the down leg. Note that the down leg does not require a duty cycle margin as large as the up leg, as only the semi-major axis is tracked.  Using this $DC^{PMDT}$, the PMDT generates a fuel optimal reference with $\Delta v = \SI{134.545}{\meter\per\second} $ and $TOF = \SI{216.957}{}$ days for the down leg. Note that for $DC^{PMDT} = 0.5$, the PMDT yields a $\Delta v = \SI{134.545}{\meter\per\second} $ and $TOF = \SI{195.261}{}$ days transfer, which has the same fuel cost but a shorter duration. However, this gives no margin of thrust to adjust for errors.  \par

\subsubsection{Down Leg Results}

\begin{figure}[hbt!]
  \centering
  \begin{subfigure}[b]{0.49\textwidth}
 \centering
    \includegraphics[width = \textwidth]{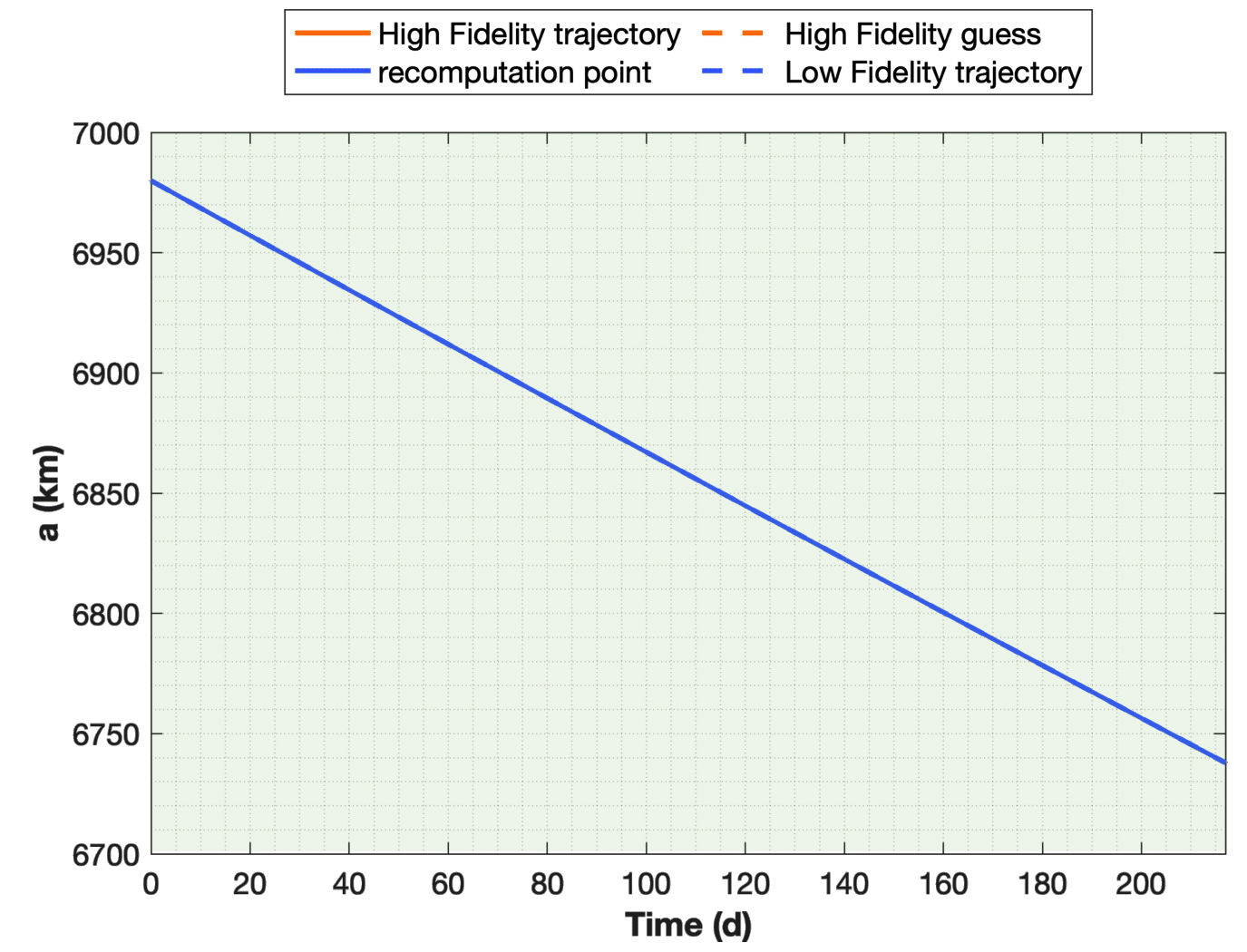}
    \caption{The fuel optimal solution with the high and low fidelity spacecraft propagations}
    \label{dlegne1}
  \end{subfigure}
  \hfill
  \begin{subfigure}[b]{0.49\textwidth}
   \centering
    \includegraphics[width = \textwidth]{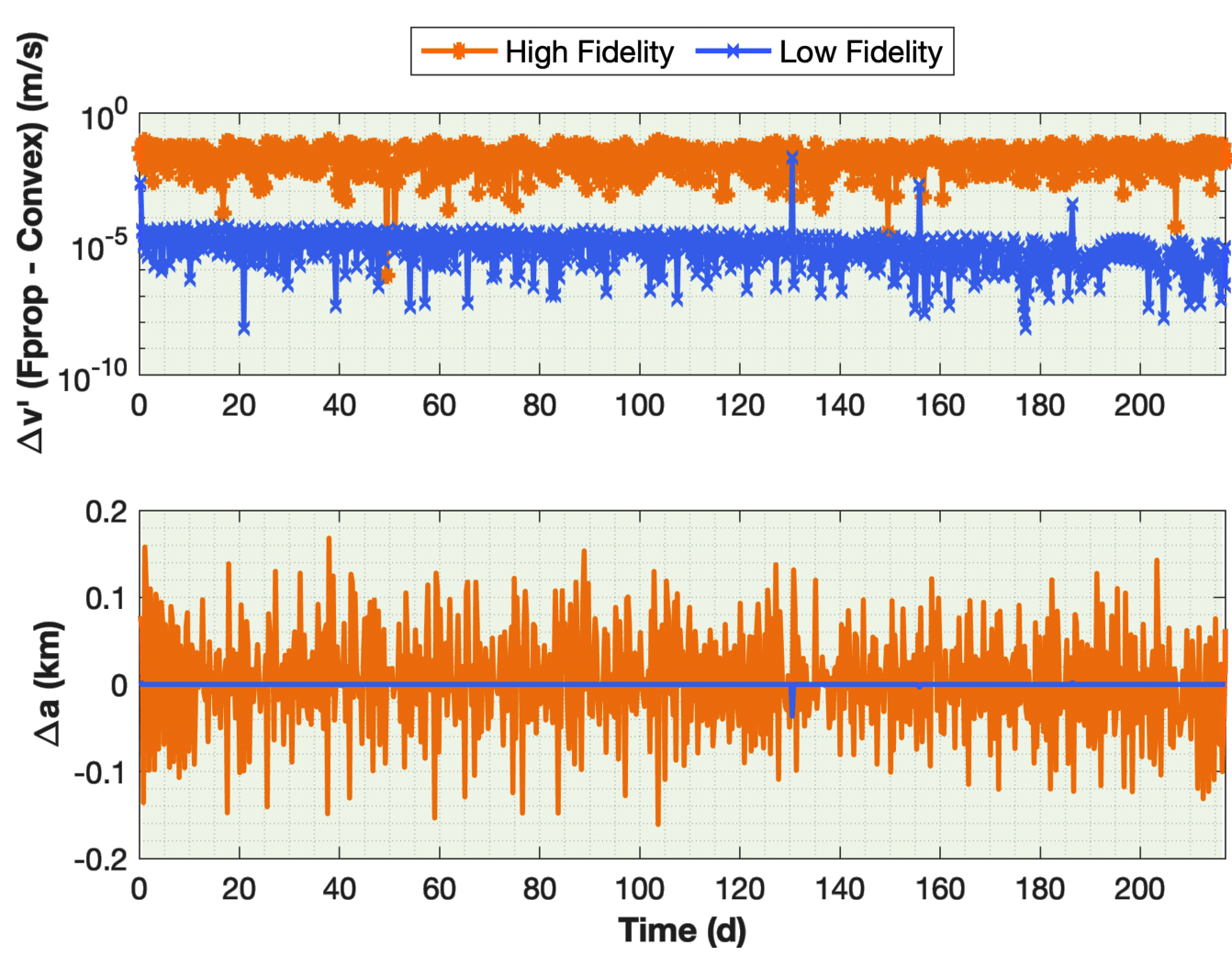}
    \caption{The evolution of $\Delta v'$ and  semi-major axis errors with respect to the reference trajectories over time}
    \label{dlegne2}
  \end{subfigure}
\caption{Down leg fuel-optimal results (perfect thrust)}
  \label{fig:main4}
\end{figure}

Figures \ref{dlegne1} and \ref{dlegne2} show the results with perfect thrust, which illustrates that good semi-major axis tracking accuracy can be maintained throughout the trajectory with both high and low fidelity dynamics. The low-fidelity dynamics experience less deviation from the reference than the high-fidelity dynamics, as expected. \par 

\begin{figure}[hbt!]
    \centering
    \includegraphics[width = \textwidth]{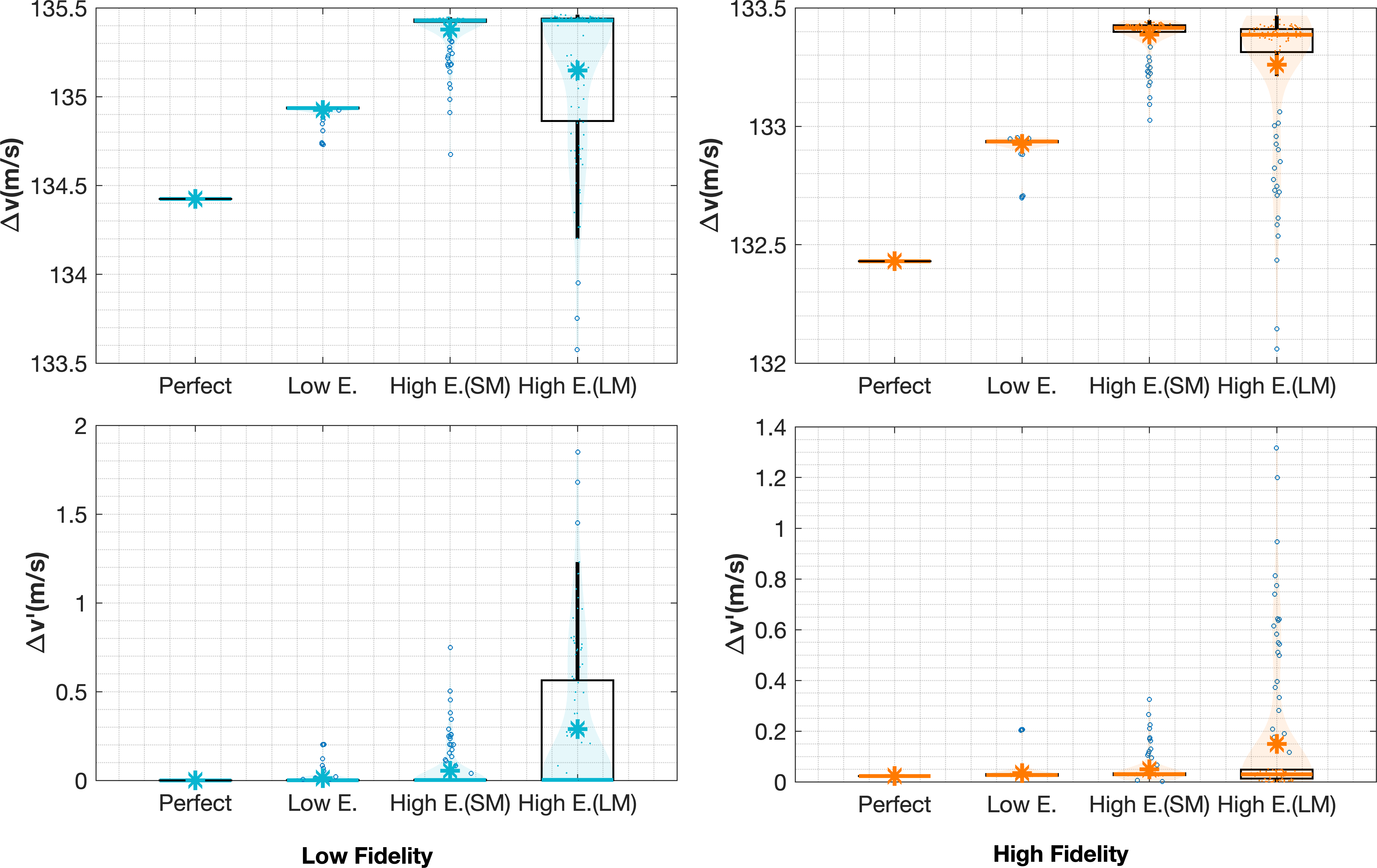}
    \caption{Down leg fuel-optimal results with thrust errors (SM: short misthrust, LM: long misthrust)}
    \label{dlegh}
\end{figure}

Once again, Monte-Carlo simulations of 100 samples were conducted for each of the low and high error settings. These results are given in Figure \ref{dlegh}.  The numerical results of the perfect thrust case and the  Monte-Carlo simulations ($25\%$ and $75\%$ quantiles) are given in Table \ref{td}. 
\par
Similar to the up leg, due to the discrepancy of dynamics between the convex tracking and forward propagation of the spacecraft, the $\Delta v'_{fp}$s of the high-fidelity simulations are greater than those of the low-fidelity simulations. From Table \ref{td}, it can be seen that the fuel consumptions ($\Delta v$s) of the high-fidelity simulations are lower than that of the low-fidelity simulations for all error cases, unlike the up leg. This is due to the discrepancy of the drag models used.  As mentioned earlier, the NRLMSISE-00 model used in the high fidelity dynamics estimates higher densities than the HP model used in the PMDT and the low fidelity dynamics, especially for low altitudes \cite{compare}. Hence, unlike in the up leg, the higher drag of the high-fidelity dynamics has a positive net effect on the down trajectory, reducing the $\Delta v$ required to lower the spacecraft's altitude.   Note that while some fuel must be used to counteract the additional perturbations of the high-fidelity dynamics, the opposing impact of the drag discrepancy is evidently much stronger.\par

Also, unlike the up leg, a long misthrust event is unlikely to cause a recomputation in the down leg, as only the semi-major axis is tracked. Due to this reason and the absence of phase angle matching, for all cases in the Monte-Carlo simulations of the down leg, the time of flight remains the same as the PMDT guess. However, when misthrust events happen towards the end of the transfer, there is insufficient time to perfectly correct the trajectory to match up with the reference, resulting in the low $\Delta v$ and high $\Delta v'_{fp}$  {asymmetry} seen in Figure \ref{dlegh}.\par 

Lastly, Figure \ref{dlegh} shows that the $\Delta v$ and $\Delta v'_{fp}$ distributions in the high-fidelity, high-error, long misthrust simulations are significantly less skewed compared to their low-fidelity, high-error counterparts. This can be attributed to the additional reduction of altitude in the high-fidelity dynamics during long misthrust events due to the greater atmospheric drag resultant of the NRKMSISE-00 model. Consequently, even when a long misthrust event occurs towards the later stages of the transfer, the high-fidelity simulation is capable of a more rapid recovery without the need for additional fuel. \par 

\begin{table}[hbt!]
\centering
\caption{Down leg fuel-optimal results \protect\footnotemark}
\label{td} 
\adjustbox{max width=\textwidth}{%
\begin{tabular}{cccc} \hline
Spacecraft Propagation & $\Delta \bar{a}$ (km)  & $\Delta v$ (m/s) & $\Delta v^{'}_{fp}$ (m/s) \\ \hline 
\multicolumn{4}{c}{PMDT}                                                                                                      \\
- & - & 134.56& -\\ \hline 
\multicolumn{4}{c}{Perfect thrust  (Single iteration)}                                                                                                      \\
 Low Fidelity & 3.1e-05 & 134.42 & 1.77e-05 \\
 High Fidelity & -0.0398 & 132.43  & -0.0228  \\\hline 
\multicolumn{4}{c}{Low Error [$Q_1 \ Q_3$]}                                                                                                      \\
Low Fidelity & [0.000335 0.00354] & [134.93 134.94] & [0.000558 0.00203] \\
High Fidelity & [-0.0513 -0.0424]& [132.93 132.94] & [0.0244 0.0308] \\ \hline  
\multicolumn{4}{c}{High error (short misthrust) [$Q_1 \ Q_3$]}                                                                                                      \\
Low Fidelity & [0.000747 0.00629]& [135.42 135.43] & [0.00101 0.00357] \\
High Fidelity & [-0.0574 -0.0438] & [133.40 133.43] & [0.0265 0.034] \\ \hline 
\multicolumn{4}{c}{High error (long misthrust) [$Q_1 \ Q_3$]}   \\
Low Fidelity & [0.00174 0.989]& [134.86 135.44] & [0.00138 0.562] \\
High Fidelity & [-0.0518 0.0656]&  [133.31 133.41] & [0.0131 0.0481] \\\hline
\end{tabular}}

\end{table}

\footnotetext{Note that as no phase tracking is conducted in the down leg, and as no recomputations are required in any of the cases with errors, the $TOF$ always remains at the PMDT output of $216.957$ days.}

\section{Conclusion}
Providing autonomous guidance for ADR missions poses a significant technical obstacle as it requires performing complex optimizations in real-time while also accounting for unexpected situations such as thruster misalignments and misthrust events. This paper proposes a convex optimization-based MPC approach to address these difficulties. In this algorithm, the reference trajectory is updated in an MPC manner using the highly accurate PMDT method, while tracking is conducted using convex optimization. The proposed guidance is shown to provide more accurate control compared to classical guidance laws as it accounts for perturbations in real time and adapts to significant divergences. It can also conduct phase matching alongside the transfer when required, such that a target phasing angle can be met at the end of a transfer. The MPC guidance can recompute a new optimal trajectory upon encountering large deviations from the current reference, which is shown to help resolve the accumulation of modeling and convexification errors. The results from the Monte-Carlo guidance simulations reveal that the spacecraft can closely adhere to the optimized reference even in the presence of thrust errors and misthrust events. 
This holds true even when forward propagation involves more complex, high-fidelity dynamics. The comparison between single iteration and successive convex tracking shows that successive convexifications significantly increase computational time but only provide a marginal increment in accuracy, underscoring the rationale for not utilizing successive convexification in the guidance algorithm. The findings also show that the spacecraft can recompute a new optimal trajectory and still successfully reach its target when deviated significantly from the original reference. 

\appendix
\section*{Appendix}

\subsection{Preliminary Mission Design Tool (PMDT) Algorithms}
The RAAN matching algorithm, the Extended Edelbaum algorithm and the PMDT forward propagation algorithm from \cite{ADRMW} are given here for reference. 
\begin{algorithm}[hbt!]
\caption{Extended Edelbaum method \cite{ADRMW}}\label{A1}
\textbf{Input} Initial and final mean semi-major axes $\bar{a}_0,\bar{a}_f$, initial and final mean inclinations $\bar{i}_0 , \bar{i}_f$, max. thrust $T_{max}$. duty cycle $DC$, initial spacecraft mass $m_0$
\begin{algorithmic}
\State Calculate the initial and final orbital velocities 
\begin{align} 
  \bar{V}_0 = \sqrt{\mu/\bar{a}_0},\bar{V}_f = \sqrt{\mu/\bar{a}_f} 
\end{align}
\State Calculate thrust acceleration $\mathscr{a}_{t}  = {T_{max}}/{m_0}$.
\State Calculate $\Delta v$ and mission time of flight ($TOF$) using Eqs. \eqref{eq1} and \eqref{eq2}.
\begin{align} 
  \Delta v_{\text{total}} &= \sqrt{\bar{V}_0^2 + \bar{V}_f^2 - 2\bar{V}_0 \bar{V}_f \cos(\pi/2 \Delta \bar{i})} \ \text{where} \  \Delta \bar{i} =\bar{i}_f - \bar{i}_0  \label{eq1}\\ 
    TOF &= { \Delta v_{\text{total}}}/\mathscr{a}_{t}  \label{eq2}
\end{align}
\State Calculate initial yaw steering angle $\beta_0$ using Eq. \eqref{eq3}. 
\begin{equation}\label{eq3}
\tan{\beta_0} = {\sin{(\pi/2 \Delta \bar{i} )}}/\left(\frac{\bar{V}_0}{\bar{V}_f} - \cos(\pi/2 \Delta \bar{i}) \right)
\end{equation}
\State Discretize $TOF$ into N segments and compute semi-major axis, inclination, yaw steering angle, and $\Delta v$ per segment using Eqs. \eqref{eq4}, \eqref{eq5}, \eqref{eqbeta}, and \eqref{eq2}.
\begin{align} 
\bar{a}(t) &= {\mu}/({\bar{V}_0^2 + \mathscr{a}_{t}^2  t^2 - 2 \bar{V}_0 \mathscr{a}_{t}  t\cos(\beta_0)})
 \label{eq4}\\ 
\bar{i}(t ) &= \bar{i}_0 +\textrm{sgn}{(\bar{i}_f-\bar{i}_0)}\frac{2}{\pi}\left[\tan^{-1}\left(\frac{\mathscr{a}_{t}  t - \bar{V}_0 \cos{\beta_0}}{\bar{V}_0 \sin{\beta_0}}\right) + \frac{\pi}{2}-\beta_0 \right] \label{eq5}\\
\beta (t ) &= \tan^{-1}({\bar{V}_0 \sin{\beta_0}}/({\bar{V}_0\cos{\beta_0} -\mathscr{a}_{t}  t})) \label{eqbeta}
\end{align}
\For{$k = 1:N$}
    \State  Calculate mass $m(k)$ using $\Delta v(k)$ via rocket equation.
    \State Calculate $\mathscr{a}_{t}^{PMDT} (k) = T_{max}/m(k)$
    \State Compute new transfer time using Eq. \eqref{eq6}. 
    \begin{equation}
    t (k+1) = t (k) + \frac{\Delta v(k+1)-\Delta v(k)}{\mathscr{a}^{PMDT}_t (k) DC}
    \label{eq6}
    \end{equation}
    \State Calculate drag acceleration ($\mathscr{a}_{drag}$) at $t_k$ using Eq. \eqref{eq7}. 
    \begin{equation}\label{eq7} 
    \mathscr{a}_{drag}(k) =-0.5 {\rho C_d A v(k)^2}/{m(k)}    
    \end{equation}
    \State where $\rho, C_d, A, v(k)$, and $m(k)$ represent the air density, drag coefficient, frontal area, velocity, and mass. 
    \State {Calculate semi major axes $a_d({k+1})$ from $ \mathscr{a}_{drag}(k)$ using Eq. \eqref{eq4}, assuming $ \mathscr{a}_{drag}(k)$ remains constant from $t_{k}$ to $t_{k+1}$.}  
    \State {Update $\bar{a}({k+1}) \leftarrow \bar{a}({k+1}) - (\bar{a}(k) - \bar{a}_d(k+1) $).} 
    \State Return to the first step of this algorithm and repeat the procedure from $t({k+1})$ to $t_{f}$. 
    \State Propagate the RAAN using Eqs. \eqref{RAAN1} and \eqref{RAAN2}
    \begin{align} 
  \dot{\bar{\Omega}}(k)&=-\frac{3}{2} J_2 \sqrt{\frac{\mu}{\bar{a}(k)^3}} \left(\frac{R_e}{\bar{a}(k)} \right)^2 \cos{\bar{i}} \label{RAAN1}\\ 
    \bar{\Omega}{({k+1})} &=  \bar{\Omega}{({k})} +  \dot{\bar{\Omega}}(k)(t({k+1})-t({k})) \label{RAAN2}
\end{align}
\EndFor  
\end{algorithmic}
\textbf{Output}  $TOF$, $\Delta v$, $t, \bar{a}(t), \bar{i}(t)$, $\bar{\Omega}(t)$  and ${\mathscr{a}_{t}^{PMDT}}$  
\end{algorithm}

\begin{algorithm}[hbt!]
\caption{RAAN matching method \cite{ADRMW}}\label{A2}
\textbf{Input} Initial orbit elements ($\bar{a}_0$, $\bar{i}_0$, $\bar{\Omega}_{t_0,0}$), target orbit elements ($\bar{a}_T$, $\bar{i}_T$, $\bar{\Omega}_{t_0,T}$), and drift orbit elements ($\bar{a}_d, \bar{i}_d$).
\begin{algorithmic}
\State \textbf{Thrust phase 1}
\State{{Input: $\bar{a}_0$, $\bar{i}_0$, $\bar{a}_d$, and $\bar{i}_d$}}
\State Calculate $TOF_{T1}$, $\Delta v_{T1}$, $t_{T1}, \bar{a}(t_{T1}), \bar{i}(t_{T1}), \bar{\Omega}(t_{T1})$ and ${\mathscr{a}_{t,T1}^{PMDT}}$  using Algorithm \ref{A1}. 
\State  Calculate $\Delta \bar{\Omega}_{T1} = \bar{\Omega}(t_{T1}) -\bar{\Omega}_{t_0, 0} $  (RAAN change of the spacecraft due to precession during thrust phase 1.)   
\State Output: 
\State  \textbf{Thrust phase 2}
\State{{Input: $\bar{a}_d$, $\bar{i}_d$, $\bar{a}_T$, and $\bar{i}_T$}}
\State Calculate $TOF_{T2}$, $\Delta v_{T2}$, $t_{T2}, \bar{a}(t_{T2}), \bar{i}(t_{T2}), \bar{\Omega}(t_{T2})$ and ${\mathscr{a}_{t,T2}^{PMDT}}$  using Algorithm \ref{A1}. 
\State  Calculate $\Delta \bar{\Omega}_{T2} = \bar{\Omega}(t_{T2}) -\bar{\Omega}_{t_0, T} $  (RAAN change of the spacecraft due to precession during thrust phase 2.)  
\State \textbf{Drifting}
\State{{Input: $\bar{a}_d, \bar{i}_d , \bar{a}_T, \bar{i}_T,TOF_{T1},TOF_{T2}, \Delta \bar{\Omega}_{T1}, \Delta \bar{\Omega}_{T2}, \bar{\Omega}_{t_0,0}$ and $\bar{\Omega}_{t_0,T}$}}
\State Calculate the drift rate of the spacecraft $\dot{\bar{\Omega}}_{s/c}$ and the drift rate of the target $\dot{\bar{\Omega}}_{T}$ using Eq. \eqref{RAAN1}.
\State Calculate $TOF_d$ (Drifting time required to match with the final RAAN) using Eq. \eqref{tw}, which equates the RAAN reached by the Servicer to the RAAN of the debris at arrival time.
{
\begin{equation}\label{tw}
    \bar{\Omega}_{t_0,0}+ \Delta \bar{\Omega}_{T1} +\Delta \bar{\Omega}_{T2} + \dot{\bar{\Omega}}_{s/c} TOF_d = \bar{\Omega}_{t_0, T}+ \dot{\bar{\Omega}}_{T}(TOF_d+TOF_{T1}+TOF_{T2}) + 2k\pi
\end{equation}}
where $k$ is an integer.
\State Calculate the RAAN reached at the end of the drift time  $\bar{\Omega}_d =   \bar{\Omega}_{t_0,0}+ \Delta \bar{\Omega}_{T1}+ \dot{\bar{\Omega}}_{s/c} TOF_d$ 
\State Calculate the $\Delta v$ used to offset the drag in the drifting phase ($\Delta V_d$). This is achieved by setting the thrust magnitude equal to the drag acceleration (Eq. \eqref{eq7}) acting in the opposite direction during drifting. Subsequently,  
{
\begin{equation}
    \Delta V_d = -\int^{TOF_d}_0 \mathscr{a}_{drag} \ dt 
\end{equation}}
\State{{Output: $TOF_d$ and $\Delta v_{p}$}}\par 
\State \textbf{Return}: Calculate the total $\Delta v$ and $TOF$.
\begin{align}
        \Delta v^{PMDT} &=\Delta v_{T1}+ \Delta v_{T2} + \Delta V_d \\ 
        TOF^{PMDT} &=TOF_{T1}+TOF_d+TOF_{T2} 
\end{align}
\State Optimize $\Delta v^{PMDT}$ or $TOF^{PMDT}$ by selecting the best $V_d$ and $i_d$ using the interior point method.
\end{algorithmic}
\textbf{Output} Optimized $\Delta v^{PMDT}$, $TOF^{PMDT}$, along with 
\begin{align}
    \mathscr{a}_t^{PMDT} &= [{\mathscr{a}_{t,T1}^{PMDT}}, -\mathscr{a}_{drag}, {\mathscr{a}_{t,T2}^{PMDT}}]^T \\
    t^{PMDT} &= [t_{T1} , TOF_{T1} + TOF_d, TOF_{T1} + TOF_d+ t_{T2}]^T \\
    \bar{a}^{PMDT} &= [\bar{a}(t_{T1}), \bar{a}_d, \bar{a}(t_{T2})]^T \\
    \bar{i}^{PMDT} &= [\bar{i}(t_{T1}), \bar{i}_d, \bar{i}(t_{T2})]^T \\
    \bar{\Omega}^{PMDT} &= [\bar{\Omega}(t_{T1}), \bar{\Omega}_d, \bar{\Omega}(t_{T2})]^T 
\end{align}
corresponding to the optimal transfer. 

\end{algorithm}

\begin{algorithm}[hbt!]
\caption{Forward propagation of the PMDT}\label{alg:2}
\textbf{Input} PMDT output $t^{PMDT}$ and $\beta^{PMDT}$, Initial state $[\bm{x}_{CC}(t_0),m_0]$
\begin{algorithmic}
\State Integrate $[\dot{\boldsymbol{x}}_{CC}, \dot{m}] = propagate(t,\bm{x}_{CC},m)$ in $t^{convex}$ to obtain $[\boldsymbol{x},m]$.  
\State     \hspace{\algorithmicindent}  $[\dot{\boldsymbol{x}}_{f,CC}, \dot{m}] = propagate(t,\bm{x}_{CC},m)$:
\State  \hspace{\algorithmicindent}   \hspace{\algorithmicindent}   Calculate the mean COE coordinates  $\bm{\bar{x}}_{COE} = oscCart2meanKep(\boldsymbol{x}_{CC})$.
\State   \hspace{\algorithmicindent}   \hspace{\algorithmicindent}   Calculate the mean argument of latitude: $\bar{L} ^{PMDT}(t) =  \bm{\bar{x}}_{COE}(5) + \bm{\bar{x}}_{COE}(6) $
\State  \hspace{\algorithmicindent}   \hspace{\algorithmicindent}  Calculate the eclipse/DC profile as given in \cite{ADRMW}.
\State \hspace{\algorithmicindent}   \hspace{\algorithmicindent}  \hspace{\algorithmicindent}  Calculate $q_1 = \cos^{-1}(\cos(\bar{L} ^{PMDT} - \bar{L}_c)) $  and $q_2 =  \cos^{-1}(\cos(\bar{L} ^{PMDT}- \bar{L}_c -\pi )) $ \par \Comment{$\bar{L}_c$ is the centre of the eclipse, calculated as discussed in \cite{ADRMW}}
 \State \hspace{\algorithmicindent}    \hspace{\algorithmicindent}  \hspace{\algorithmicindent}  \textbf{if} { $q_1 < \frac{\pi}{2}(1-DC^{PMDT})$ \textbf{or}  $q_2 < \frac{\pi}{2}(1-DC^{PMDT})$}  \textbf{then} $\eta=0$ 
\State  \hspace{\algorithmicindent}   \hspace{\algorithmicindent}  \hspace{\algorithmicindent}  \textbf{else} $\eta=1$ \Comment{Note that $\eta=0$ indicates drifting.}
\State \hspace{\algorithmicindent}  \hspace{\algorithmicindent}   Compute out-of-plane thrust angle $\beta$  using Eq. 6 in \cite{ADRMW}.
\State \hspace{\algorithmicindent}  \hspace{\algorithmicindent}   Reverse thrust direction across the line of nodes to make the inclination change possible.
 \State  \hspace{\algorithmicindent}   \hspace{\algorithmicindent}  \hspace{\algorithmicindent}  \textbf{if}  { $\bar{i}_f > \bar{i}_0$}
 \State  \hspace{\algorithmicindent}   \hspace{\algorithmicindent} \hspace{\algorithmicindent}  \hspace{\algorithmicindent}   \textbf{if}  { $\bar{L} ^{PMDT}(t)   > \pi/2$ or $\bar{L} ^{PMDT}(t)  \leq 3\pi/2$} \textbf{then} $\beta = -|\beta| $  \textbf{else}  { $\beta = |\beta|$}
 \State  \hspace{\algorithmicindent}   \hspace{\algorithmicindent}  \hspace{\algorithmicindent}  \textbf{else}
 \State  \hspace{\algorithmicindent}   \hspace{\algorithmicindent} \hspace{\algorithmicindent}  \hspace{\algorithmicindent}   \textbf{if} { $\bar{L} ^{PMDT}(t)  \geq \pi/2$ or $\bar{L} ^{PMDT}(t)  < 3\pi/2$} \textbf{then}
$\beta = |\beta|$  \textbf{else} 
$\beta = -|\beta|$
\State \hspace{\algorithmicindent}  \hspace{\algorithmicindent} Compute thrust: $\bm{\mathscr{a}}_{t} = \eta T_{max}/m[0, \cos{\beta}, \sin{\beta}]$, and convert to TOD. 
\State \hspace{\algorithmicindent}  \hspace{\algorithmicindent}  Compute the gravity $\bm{\mathscr{a}}_g$, $J_2$ $\boldsymbol{\mathscr{a}}_{J2}$ and drag $\boldsymbol{\mathscr{a}}_{drag}$ accelerations as in \cite{8}. 
\State \hspace{\algorithmicindent}  \hspace{\algorithmicindent}  Calculate the state derivative. $\dot{\boldsymbol{x}}  =   \boldsymbol{\mathscr{a}}_{t}+ \bm{\mathscr{a}}_g + \boldsymbol{\mathscr{a}}_{J2} + \boldsymbol{\mathscr{a}}_{drag} $.
\State \hspace{\algorithmicindent}  \hspace{\algorithmicindent}  Determine the mass derivative $\dot{m} = {\eta T_{max}}/({I_{sp} g_0})$
\end{algorithmic}
\textbf{Output}  $\bm{x}_{f,CC}$ (the final state reached at $t_f$) and   $\bar{L}^{PMDT}$ (the mean argument of latitude profile)
\end{algorithm}

\bibliographystyle{unsrt}

\end{document}